\documentstyle[aps,pra,multirow,hhline,epsfig]{revtex}      

\begin{document}     
\title{Discrete singular convolution and 
         its application to  computational
          electromagnetics}
    
\author{G.~W.~Wei\footnote{
        Tel: (65) 874 6589; Fax: (65) 774 6756; 
        E-mail: cscweigw@nus.edu.sg}}    

\address{Department of Computational Science, 
      National University of Singapore\\
      %3 Science Drive 2\\
      Singapore 117543, R. Singapore}

\date{\today}    
 
\maketitle    
    
\begin{abstract}    
A new computational algorithm, the discrete singular convolution 
(DSC), is introduced for computational electromagnetics. The basic 
philosophy behind the DSC algorithm for the approximation of functions 
and their derivatives is studied. Approximations to the delta 
distribution are constructed as either bandlimited reproducing kernels 
or approximate reproducing kernels. A systematic procedure 
is proposed to handle a number of boundary conditions which 
occur in practical applications. The unified features of the DSC 
algorithm for solving differential equations are explored from the 
point of view of the method of weighted residuals. It is demonstrated 
that different methods of implementation for the present algorithm, 
such as global, local, Galerkin, collocation, and finite difference, 
can be deduced from a single starting point. Both the computational 
bandwidth and the accuracy of the DSC algorithm are shown to be 
controllable. Three example problems are employed to illustrate the 
usefulness, test the accuracy and explore the limitation of the DSC 
algorithm. A Galerkin-induced collocation approach is used for a 
waveguide analysis in both regular and irregular domains and for  
electrostatic field estimation via potential functions. Electromagnetic 
wave propagation in three spatial dimensions is integrated by using a 
generalized finite difference approach, which becomes a global-finite 
difference scheme at certain limit of DSC parameters. Numerical 
experiments indicate that the proposed algorithm is a promising 
approach for solving problems in electromagnetics.
 

\end{abstract}

\newpage
%\hfil    
%\vspace*{2cm}
 
\section{Introduction}

     Recently, computational electromagnetics (CEM) has emerged as
a distinct  scientific  discipline for the study and understanding of 
a wide variety of  electrical and electronic engineerings 
problems\cite{BahLev,Wang,TSB,AndVol,GarSal,SLC,Honma,TsuSeb,UmaTaf,Nathan}.
As a natural extension to the analytical approach to the 
Maxwell equation, the CEM is based on numerically solving 
the governing equations in either the partial differential form  
or in the integral equation form.
The complexities of physics and of the field geometry
are no longer the limiting factors to CEM as they are to the 
analytical approach. With the advent of high-performance 
digital computers, CEM is emerging as a powerful approach for 
solving practical problems in electromagnetics. 
In fact, the dramatic progress in solving the Maxwell equation 
made in the last few decades\cite{Honma,TsuSeb,UmaTaf,Nathan} 
has opened up a new research
frontier in electromagnetics, plasmadynamics, optic engineering, 
and  an interphase between electrodynamics and quantum dynamics. 
As a computational discipline, the success of the CEM is vitally 
dependent
on the virtues of computational algorithms, such as numerical accuracy, 
stability and efficiency. These in turn depend on grid methods and 
numerical schemes for the solution of the Maxwell equation. 

A variety of computational techniques have been used for CEM,
including wavelet analysis\cite{BahLev,Wang},
integer lattice gas automata\cite{TSB},
hierarchical tangential vector finite elements\cite{AndVol},
 Nedelec tetrahedral element method\cite{GarSal} and other 
approaches\cite{SLC,UmaTaf}.
Typically, grid methods used in CEM are  either   
global\cite{Lanczos,CooTuk,FinScr,Orszag,Fornberg}, such as fast 
Fourier transform, spectral methods and pseudospectral methods, or 
local\cite{ForWas,IsaKel,Zienkiewicz,DesAbe,Oden,Nath,Fenner,Cheung,Rao,Reddy},
such as finite difference, finite volume and finite element methods.
Global methods are highly accurate but are cumbersome  
to implement in complex geometries and non-conventional 
boundary conditions. For example, a global method may  
converge slowly in a waveguide mode analysis due to 
irregular boundary conditions.
In contrast, local methods are easy to implement for  
complex geometries and discontinuous boundary conditions. 
However, the accuracy of local methods is usually very low.   
There exists many problems in CEM which require both high 
computational accuracy and high numerical flexibility in handling 
complex geometries. These problems are characterized by a geometry 
which has a large domain size, i.e., the dimensions of the scatterer
greatly exceed the wavelength of the incident wave. A typical  
example is the radar cross-section analysis of an entire 
airplane with an incident  electromagnetic wave having a frequency of 
the order of ten GHz. To deal with such problems, 
it is desired to have a computational method 
that has both global methods' accuracy and local methods' 
flexibility.

More recently, discrete singular convolution (DSC)     
algorithm was proposed as a potential approach      
for the computer realization of singular      
convolutions\cite{weijcp99,weijpa20}.  
Sequences of  approximations to the  
singular kernels of Hilbert type, 
Abel type and delta type were constructed.    
Applications were discussed to analytical     
signal processing, Radon transform and surface    
interpolation. The mathematical foundation      
of the DSC algorithm is the theory of     
distributions\cite{Schwartz} and the theory of wavelets.    
Numerical solutions to differential equations      
are formulated via the singular kernels of the delta type. 
By appropriately selecting parameters of a DSC kernel, 
the  DSC approach exhibits controllable accuracy      
for integration and shows excellent flexibility      
in handling complex geometries and boundary conditions. 
Many DSC kernels, such as (regularized)    
Shannon's delta kernel, (regularized) Dirichlet kernel,     
(regularized) Lagrange kernel and     
(regularized) de la Vall\'{e}e      
Poussin kernel, have been constructed\cite{weijcp99}. 
Practical applications were examined for the numerical solution of      
the Fokker-Planck equation\cite{weijcp99,weijpa20} and  for   
the Schr\"{o}dinger equation\cite{weijpb2000}.   
Another development in the application of the     
DSC algorithm is its use in computing numerical solutions of   
the Navier-Stokes equation
and in structural analysis\cite{weicmame20}.      
In the context of image processing, DSC kernels were  
used to facilitate a new anisotropic diffusion operator  
for image restoration from noise\cite{weiieee}.  
Most recently, the DSC algorithm was used to   integrate the      
(nonlinear) sine-Gordon equation      
with the initial values close to a homoclinic      
manifold singularity\cite{weiphysica99}, for which      
conventional local methods     
encounter great difficulties and    
result in numerically induced chaos\cite{Ablowitz}.

The purpose of this paper is to study the computational 
philosophy of the DSC algorithm and to introduce
the algorithm for computational electromagnetics (CEM).
For the purpose of numerical computation, both bandlimited 
reproducing kernels and approximate reproducing kernels are 
discussed as sequences of approximations to the universal 
reproducing kernel, the delta distribution.
A systematic treatment is proposed for 
handling a general class of boundary conditions. 
We explore the unified feature of the DSC algorithm
for numerical approximation of differential 
equations. It is found that several conventional 
computational methods, such as  methods of global, 
local, Galerkin, collocation, and  finite difference
can be derived from a single starting point.
In particular, a Galerkin-induced collocation 
algorithm is discussed. A set of generalized finite
difference schemes are shown to exhibit global-finite
difference features at certain limit of DSC parameters.
The present algorithm is shown to have controllable accuracy.
The potential of the DSC algorithm for computational 
electromagnetics is explored by using three classes of 
problems, the eigenmode analysis of waveguide, the potential 
function analysis of electrostatics and the propagation of 
electromagnetic waves.

This paper is organized as follows:
In Section II, we study the computational 
philosophy of the DSC algorithm. 
A number of new DSC kernels are constructed
as approximations to the universal reproducing 
kernel---the delta distribution.
Approximation of functions and their derivatives 
is discussed. A systematic treatment of boundary 
conditions is proposed for being used in 
implicit schemes. The capability 
of the DSC algorithm is analyzed for solving 
differential equations in Section III. The unified 
feature of the DSC algorithm is explored in the 
framework of the method of weighted residuals.
The application of the DSC algorithm to CEM
is introduced in Section IV. 
The utility and robustness of the proposed method 
is illustrated by  a few numerical experiments.
This paper ends with a conclusion.

\section{The Discrete Singular Convolution}

\subsection{Approximation of singular convolution}

Singular convolutions (SC) are a special    
class of mathematical transformations, which appear
in many science and engineering  problems,    
such as  Hilbert transform, Abel transform and Radon transform.  
It is most convenient to discuss the singular convolution 
in the context of the theory of    distributions. 
The latter has a significant impact in mathematical   
analysis. Not only it provides a rigorous justification    
for a number of informal manipulations in physical and    
engineering sciences, but also it opens a new area    
of mathematics, which in turn gives impetus to many other    
mathematical disciplines, such as operator calculus,    
differential equations, functional analysis,  
harmonic analysis and transformation theory.  
In fact, the theory of wavelets and frames, a new 
mathematical branch developed in recent years, 
can also find its root in the theory of distributions.

Let $T$ be a distribution and $\eta(t)$ be an
element of the space of test functions. 
A singular convolution is defined as
\begin{equation}\label{e1}
F(t)=(T\ast \eta)(t)
=\int_{-\infty}^{\infty}T(t-x)\eta(x)dx.
\end{equation} 
Here $T(t-x)$ is a singular kernel. 
Depending on the form of the kernel $T$, the 
singular convolution is the central issue for a wide 
range of science and engineering problems. 
For example, singular kernels of 
the  Hilbert type have a general form of 
\begin{equation}\label{t1}
T(x)={1\over x^n},~~~~(n=1,2,\cdots).
\end{equation}
Here, kernel $T(x)={1\over x}$ commonly occurs in 
electrodynamics, theory of linear response, 
signal processing, theory of analytic functions,  
and the Hilbert transform. When $n=2$, 
$T(x)={1\over x^2}$ is the kernel used in tomography.
Another interesting example is singular 
kernels of the Abel type
\begin{equation}\label{t1-2}
T(x)={1\over x^\beta},~~~~(0<\beta<1).
\end{equation}
These kernels can be recognized as the special cases
of the singular integral equations of Volterra type of 
the first kind. Singular kernels of the Abel type 
have applications in the area of 
holography and interferometry with phase objects 
(of practical importance in aerodynamics, heat and mass
transfer, and plasma diagnostics).
They are intimately connected with the 
Radon transform, for example, in determining the refractive 
index from the knowledge of a holographic interferogram. 
The other important example is  singular kernels 
of the  delta type
\begin{equation}\label{t2}
T(x)=\delta^{(n)}(x), ~~~ (n=0,1,2,\cdots).
\end{equation}
Here, kernel $T(x)=\delta(x)$ is of particular 
importance for 
interpolation of surfaces and curves 
(including atomic, molecular and biological 
potential energy surfaces, engineering surfaces
and a variety of image processing and pattern 
recognition problems involving low-pass filters). 
Higher-order kernels,  
$T(x)=\delta^{(n)}(x),~~(n=1,2,\cdots)$ 
are essential for numerically solving 
partial differential equations and 
for image processing, noise estimation, etc. 
However, since these kernels are singular, they cannot 
be directly digitized in computers. Hence, the singular 
convolution, (\ref{e1}), is of little numerical merit. 
To avoid the difficulty of using singular expressions 
directly in computer, we construct sequences of 
approximations ($T_\alpha$) 
to the distribution $T$ 
\begin{equation}\label{e2}
\lim_{\alpha\rightarrow\alpha_0} T_\alpha(x)\longrightarrow T(x),
\end{equation} 
where $\alpha_0$ is a generalized limit. 
Obviously, in the case of $T(x)=\delta(x)$, each element in
the sequence, $T_\alpha(x)$, is a delta sequence kernel. 
Note that one retains the delta distribution at the 
limit of a delta sequence kernel.
Computationally, the Fourier transform of the delta distribution 
is unity. Hence, it is a {\it universal reproducing kernel}
for numerical computations and an {\it all pass filter}
for image and signal processing. Therefore, the delta distribution 
can be used as a starting point for the construction of 
either band-limited reproducing kernels or approximate 
reproducing kernels. By the Heisenberg uncertainty principle,
exact reproducing kernels have bad 
localization in the time (spatial) domain, 
whereas, approximate 
reproducing kernels can be localized in both 
time and frequency representations.   Furthermore,
with a sufficiently smooth approximation, it makes
sense to consider a {\it discrete singular convolution} (DSC)
\begin{equation}\label{e3}
F_\alpha(t)=\sum_k T_\alpha(t-x_k)f(x_k),
\end{equation} 
where $F_\alpha(t)$ is an approximation to $F(t)$
and $\{ x_k\}$ is an appropriate set of discrete points on 
which the DSC (\ref{e3}) is well defined. Note that, the 
original test function $\eta(x)$ has been replaced by 
$f(x)$. The mathematical property or requirement
of $f(x)$ is determined by the approximate kernel 
$T_{\alpha}$. In general, the convolution is 
required being Lebesgue integrable.
In the rest of this paper, the emphasis is on the 
singular kernels of the delta type, their approximation,
and numerical implementation.

\subsection{Singular kernels of delta type}

The delta distribution or so called Dirac delta 
function ($\delta$) is a generalized function which is 
integrable inside a particular interval 
but in itself need not have a value. 
Heaviside introduced both the unit step Heaviside function and
the Dirac delta function as its derivative and referred to 
the latter as the  unit impulse.  
Dirac, for the first time,  
explicitly discussed the properties of $\delta$ in his classic
text on quantum mechanics; for this reason $\delta$ is often 
called Dirac delta function. However, delta distribution has 
a history which antedates both Heaviside and Dirac. It appeared 
in explicit form as early as 1822, in Fourier's 
{\it Th\'{e}orie Analytique de la Chaleur}.    
The work of Heaviside, and subsequently of Dirac, in the systematic
but informal exploitation of the step and 
delta function has made delta distribution  
familiar to physicists and engineers before Sobolev, 
Schwartz\cite{Schwartz}, 
Korevaar\cite{Korevaar} and others put it 
into a rigorous mathematical form. 
The Dirac delta function is the most important 
special case of distributions or generalized functions.  
There are three parallel descriptions for the theory of 
distributions.  
One description of distributions  is to  
characterize them as an equivalence 
class, or as a generalized limit of various Cauchy sequences 
(fundamental sequences)  and fundamental
families as rigorously defined by Korevaar\cite{Korevaar}.
This approach is particularly convenient for the delta distribution.   
Another description is to formulate them as continuous 
linear functionals on 
the space of test functions  as introduced by
Schwartz\cite{Schwartz}. The vector space of test functions is obtained 
from a class of test functions with compatible convergence or topology. 
The other description is based on generalized derivatives of 
integrable functions. Generalized derivatives are distributions rather
than well-behaved functions. The first
description is intuitive and convenient for various applications. 
The second description is particularly
elegant and concise. It is also very convenient for 
higher dimensional applications. 
The third description is useful for certain practical applications 
involving derivatives and antiderivatives.
These three  descriptions are formally equivalent and 
are commonly used for describing not only for the delta
distribution, but also for distributions in general.
The use of  many delta sequences as probability density estimators 
was discussed by Walter and 
Blum\cite{WalBlu} and others\cite{Winter,Wahba,KroTar}.

{\it Definition 1}.
The delta distribution, or so called Dirac delta function
is given as a continuous linear functional on the space of test functions,
${\cal D}(-\infty,\infty)$, 
\begin{equation}
<\delta,\phi>=\delta(\phi) 
=\int_{-\infty}^{\infty}\delta\phi=\phi(0). 
\end{equation} 
A delta sequence kernel, $\{\delta_{\alpha}(x) \}$, is a sequence of kernel 
functions on $(-\infty,\infty)$ 
which is integrable over every compact domain and their inner product 
with every test function $\phi$ converges to
the delta distribution 
\begin{equation}\label{del22}
\lim_{\alpha \rightarrow
\alpha_0} \int_{-\infty}^{\infty}\delta_{\alpha}\phi=<\delta,\phi>, 
\end{equation} 
where the (real or complex ) parameter $\alpha$ 
approaches $\alpha_0$ which can either be $\infty$ or a limit value,
depending on the situation (such a convention for $\alpha_0$ is 
used thorough out this paper).
If $\alpha_0$ represents a limit value, 
the corresponding 
delta sequence kernel is  a fundamental family. 
Depending on the explicit form of $\delta_{\alpha}$,
the condition on $\phi$ can be relaxed.
For example, if $\delta_{\alpha}$ is given as
\begin{equation}
\delta_{\alpha}(x)=\left\{\begin{array}{ll}\alpha &\mbox{for $0<x<1/\alpha~~~~
         \alpha=1,2,\cdots $}\\
       0 &\mbox{otherwise}\end{array}\right.,
\end{equation} 
then Eq. (\ref{del22}) makes sense for every $\phi$ in 
$C(-\infty,\infty)$.

There are many delta sequence kernels 
arising in the theory of partial differential equations,
Fourier transforms and signal analysis, with completely different
mathematical properties.
It is useful to have a classification of various delta 
sequence kernels for discussion, application and for new construction.  
Delta sequence kernels of positive and Dirichlet
type have very distinct mathematical properties and can serve as 
the basis of a good classification scheme. 
In particular, there is a close relation between the 
delta sequence kernel of positive type and statistical distribution 
functions. In fact, every statistical distribution function 
can be regarded as an element of a delta sequence kernel of the 
positive type. An ordinary element of delta sequence kernel
of Dirichlet type has the well-known feature of ``small wave''. 
In other words it is readily related to the wavelet scaling function. 
Moreover, classifying delta sequence kernels according to 
Schwartz class or non-Schwartz class is also 
very useful for various physical and engineering applications. 
In particular, all physically realizable states, 
either in the sense of quantum mechanics  or  classical mechanics, 
 belong to the Schwartz class\cite{Bohm}. 
Moreover, for the purpose of numerical applications to 
ill-posed problems, delta sequences of the Schwartz class 
are applicable to a wide class of functions and distributions.   
In the following two subsections, 
Delta sequence kernels of positive type and Dirichlet
type are studied.

\subsection{Delta sequence kernels of positive type }

{\it Definition 2}. Let $\{\delta_\alpha \}$ be a 
sequence of kernel functions on $(-\infty,\infty)$
which are integrable over every bounded interval. We call 
$\{\delta_\alpha \}$ a delta sequence kernel of positive type if 
\begin{enumerate}        
\item $\int_{-a}^{a}\delta_{\alpha}\rightarrow 1$~~ {\rm as }~~
     $\alpha\rightarrow \alpha_0 $~~{\rm for some finite constant $a$}.
\item For every constant $\gamma>0$, 
 $\left(\int_{-\infty}^{-\gamma}+\int^{\infty}_{\gamma}\right)
 \delta_{\alpha}\rightarrow 0$ as $\alpha\rightarrow \alpha_0 $.
\item  $\delta_{\alpha}(x)\geq 0$ for all $x$ and $\alpha$.
\end{enumerate}

{\it Example 1.  Delta sequence kernel of impulse functions}

To approximate idealized  physical concepts such as the force density of 
a unit force at the origin $x=0$, or a unit impulse at time $x=0$,
a sequence of functions given by  
\begin{equation}
\delta_{\alpha}(x)=\left\{\begin{array}{ll}\alpha &\mbox{for $0<x<1/\alpha~~~~
         \alpha=1,2,\cdots $}\\
       0 &\mbox{otherwise}\end{array}\right.
\end{equation} 
is a DSC delta sequence kernel provided $\alpha\rightarrow\infty$.
This is a commonly used density estimator in science and engineering.

{\it Example 2. Gauss' delta sequence kernel}

In the study of heat equation, Gauss' delta sequence kernel
\begin{equation}
\delta_{\alpha}(x)={1\over \sqrt{2\pi}\alpha} e^{-x^2/2\alpha^2}~~ {\rm for}~~
       \alpha\rightarrow 0
\end{equation} 
arises naturally as a distribution solution or so called weak solution. 
Gauss' delta sequence kernel has various interesting properties with regard
 to differentiability, boundedness and Fourier transforms, and it
is used to generate the ``Mexican hat'' wavelet.

{\it Example 3. Lorentz's delta sequence kernel}
                
Lorentz's delta sequence kernel
\begin{equation}\label{Lore}
\delta_{\alpha}(x)={1\over \pi} {\alpha\over x^2+\alpha^2}~~ {\rm for}~~
       \alpha\rightarrow 0
\end{equation} 
is known for its role in representing the solution of 
Laplace equation in the upper half plane.
It is commonly seen in integral equations involving the Green's function 
of the kinetic energy operator (in the momentum representation).               
It is also the expression for the line shape of various spectroscopies
when the relaxation is an exponential one in the time domain.  
A generalized expression can be written as,
\begin{equation}
\delta_{\alpha,n}(x)={1\over \pi} 
{\alpha^nx^{n-1}\over x^{2n}+\alpha^{2n}}~~ {\rm for}~~
       \alpha\rightarrow 0, ~{\rm and}~n\geq 1.  
\end{equation} 
This includes Eq. (\ref{Lore}) as a special case.

{\it Example 4. Landau's delta sequence kernel}

In the  discussion of  convergence properties of polynomial 
approximations, 
Landau introduced a delta sequence kernel
\begin{equation}
L_{n}(x)={(a^2-x^2)^n\over\int_{-a}^a(a^2-y^2)^n dy} 
~~{\rm for}~~~ n=0,1,2,\cdots ~{\rm and }~a>0. \nonumber
\end{equation} 
It becomes a delta sequence kernel
\begin{equation}
\delta_{n}(x)=\left\{\begin{array}{ll}L_n(x) 
         &\mbox{for $\mid x\mid\leq a $}\\
       0 &\mbox{otherwise}\end{array}\right.
\end{equation} 
as $n \rightarrow\infty$.
This is called Landau's delta sequence kernel.
Wavelets generated from Landau's delta sequence kernel 
can be very useful for a sufficiently large $n$.

{\it Example 5. Poisson's delta sequence kernel family }

The function given by the summation of  an infinite  series 
\begin{eqnarray}
P_{\alpha}(x)&=& {1\over\pi}\left[{1\over2}+\alpha \cos(x)+ \alpha^2 \cos(2x)+
        \cdots\right] \nonumber\\
   &=&{1-\alpha^2\over 2\pi (1-2\alpha\cos(x)+\alpha^2)},\nonumber
\end{eqnarray} 
where $0\leq \alpha<1 $ and $(-\infty <x<\infty )$, is called the 
Poisson kernel, which plays an important role in Poisson's integral formulae. 
Poisson's delta sequence kernel family is given by
\begin{equation}
\delta_{\alpha}(x)=\left\{\begin{array}{ll}P_\alpha(x) 
         &\mbox{for $\mid x\mid\leq \pi $}\\
       0 &\mbox{otherwise}\end{array}\right.
\end{equation} 
as $\alpha\rightarrow1$.                          
The Poisson's delta kernel family has a connection with the solution of 
Laplace equation in the unit disc (i.e., Dirichlet 
problem for the unit disc).

{\it Example 6. Fej\'{e}r's delta sequence kernel}

The partial sum of the discrete Fourier series  
\begin{eqnarray}\label{Dirichlet} 
D_{k}(x)&=& {1\over\pi}\left[{1\over2}+\cos(x)+ \cos(2x)+
        \cdots+\cos(kx) \right] \nonumber\\
   &=&{\sin[(k+{1\over2})x] \over 2\pi \sin{1\over2}x}, ~~~k=0,1,2,\cdots
\end{eqnarray} 
is called a Dirichlet kernel. To improve  convergence for proving a
trigonometric approximation theorem, Fej\'{e}r introduced the following
arithmetic mean 
\begin{eqnarray}
F_{k}(x)&=& {1\over k}\left[D_0(x)+D_1(x)+\cdots+D_{k-1}(x) \right] \nonumber\\
   &=&{\sin^2({1\over2}kx) \over 2\pi k \sin^2({1\over2}x) }
   ~~~-\infty< x < \infty.
\end{eqnarray} 
Then Fej\'{e}r's delta sequence kernel is given by
\begin{equation}
\delta_{\alpha}(x)=\left\{\begin{array}{ll}F_\alpha(x) 
         &\mbox{for $\mid x\mid\leq \pi $}~~~{\rm for}~~ \alpha=0,1,2,\cdots \\
       0 &\mbox{otherwise}\end{array}\right.
\end{equation} 
as  $\alpha\rightarrow\infty$.
Fej\'{e}r's delta kernel has an important application in the 
theory of reproducing kernels. It also describes intensity pattern of 
light from a regular series of pinholes  in optical physics.

{\it Example 7. Generalized Fej\'{e}r's delta sequence kernel}

It is noted that Fej\'{e}r's method of generating delta sequence kernel 
is very general. Essentially, a family of arithmetic 
means of delta sequence kernels is still a delta sequence kernel. 
The resulting delta sequence kernel can be 
called a delta sequence kernel of
Fej\'{e}r type. For instant, in a similar treatment using Dirichlet's 
continuous delta sequence kernels (see next subsection), 
one obtains the following 
Fej\'{e}r's continuous delta sequence kernel
\begin{equation}
\delta_{\alpha}(x)={2\over \pi}{\sin^2(\alpha x)\over \alpha x^2}
~~~\forall x\in R.
\end{equation} 
Obviously, this is well-defined on the real line. This expression is 
related to the intensity of light diffracted by a uniform slit.

{\it Example 8. Delta sequence kernels generated by dilation}

Let $\rho\in L^1(R)$ be a non-negative function with $\int \rho(x)dx=1$, 
dilation of $\rho$ is given by
\begin{equation}
\rho_\alpha(x)={1\over \alpha} \rho({x\over\alpha}) ~~~(\alpha>0)
\end{equation}
leads to a delta sequence kernel, $\rho_\alpha\rightarrow\delta$,
 as $\alpha\rightarrow 0$.

Physically, $\rho$ can be regarded as a statistical distribution function.
This is a general procedure and {\it Examples 2} and {\it 3} fit 
into this structure. {\it Examples 1} and {\it 6} can be expressed 
in this form by appropriate
modifications (by replacing $\alpha$ with $\beta=1/\alpha$, and then 
letting $\beta\rightarrow0$).

\subsection{ Delta sequence kernels of Dirichlet type }

{\it Definition 3}. Let $\{\delta_\alpha \}$ be a sequence 
of functions on $(-\infty,\infty)$
which are integrable over every bounded interval. We call 
$\{\delta_\alpha \}$ a delta sequence kernel of the Dirichlet type if 
\begin{enumerate}        
\item $\int_{-a}^{a}\delta_{\alpha}\rightarrow 1$~~ {\rm as }~~
     $\alpha\rightarrow \alpha_0 $~~{\rm for some finite constant $a$}.
\item For every constant $\gamma>0$;
 $\left(\int_{-\infty}^{-\gamma}+\int^{\infty}_{\gamma}\right)
 \delta_{\alpha}\rightarrow 0$ as $\alpha\rightarrow \alpha_0 $,
\item There are positive constants $C_1$ and $C_2$ such that
  $$\mid \delta_{\alpha}(x)\mid\leq {C_1\over \mid x\mid} + C_2$$ 
 for all $x$ and $\alpha$.
\end{enumerate}

{\it Example 1. Dirichlet kernel} 

The most important example of a delta sequence kernel of Dirichlet
type is Dirichlet kernel
\begin{equation}
\delta_{\alpha}(x)=\left\{\begin{array}{ll}D_\alpha(x)~~~ 
         &\mbox{for $\mid x\mid\leq \pi $} ~~~~{\rm for}~~ 
         \alpha=0,1,2,\cdots\\
       0 &\mbox{otherwise}\end{array}\right.
\end{equation} 
where $D_\alpha$ is the Dirichlet kernel given by Eq. (\ref{Dirichlet}).
Dirichlet's delta sequence kernel plays an important 
role in approximation theory
and is the key element in trigonometric polynomial approximations.
In fact, it is an exact reproducing kernel for  bandlimited, 
periodic,  $L^2$  functions. Physically, it describes the diffraction of 
light passing through a regular series of 
pinholes in which the $k$th pinhole's contribution is proportional to 
$e^{ik}$.

{\it Example 2. Modified Dirichlet kernel} 

Sometimes there is some advantage in taking the last term in 
$D_{\alpha}$  with a factor of ${1\over2}$:
\begin{eqnarray}
D^*_{\alpha}(x)&=& D_{\alpha}-{1\over2\pi}\cos(\alpha x) \nonumber\\
   &=&{\sin(\alpha x) \over 2\pi \tan({1\over2}x)}, ~~~\alpha=0,1,2,\cdots.
\end{eqnarray} 
This is the so-called modified Dirichlet kernel.
The difference $D_\alpha-D^*_\alpha$ tends uniformly to zero on $(-\pi,\pi)$ as
$\alpha\rightarrow \infty$. They are equivalent with respect to convergence.

The expression given by 
\begin{equation}
\delta_{\alpha}(x)=\left\{\begin{array}{ll}D^*_\alpha(x)~~~~ 
         &\mbox{for $\mid x\mid\leq \pi $} 
~~~{\rm for}~~ \alpha=0,1,2,\cdots \\
       0 &\mbox{otherwise}\end{array}\right.
\end{equation} 
is a delta sequence kernel of Dirichlet type as $\alpha\rightarrow \infty$.

{\it Example 3. Lagrange kernel}

 Lagrange interpolation  formula   
\begin{equation} \label{L1-0}  
L_{M,k}(x)=  
\prod_{i=k-M,i\neq k}^{i=k+M}{x-x_i\over x_k-x_i},~~~(M\geq 1)  
\end{equation}   
is defined on an interval $(a,b)$ with 
a set of $2M+1$ ordered discrete points,
\begin{equation}
\{x_i\}_{i=k-M}^{k+M}:~ 
x_{k-M}=a< x_{k-M+1} < \cdots< x_k < \cdots < x_{k+M}=b. 
\nonumber
\end{equation} 
It converges to the delta distribution as
\begin{equation}\label{Lag}
a\rightarrow-\infty, b\rightarrow\infty
~~{\rm and}~~ 
\sup_{\forall x_i,x_j\in (a,b)}\mid
x_i-x_j\mid\rightarrow0.
\end{equation}             
Obviously, these limits  imply 
$M\rightarrow\infty$. 
Since the delta distribution has only a point support,
the Lagrange interpolation formula is a delta sequence 
\begin{equation}\label{Lag-dsc}
\delta_{M,k}(x)=\left\{\begin{array}{ll}L_{M,k} (x)~~~~ 
         &\mbox{for $a\leq x \leq b $} 
~~~{\rm for}~~ M=1,2,\cdots \\
       0 &\mbox{otherwise}\end{array}\right.
\end{equation} 
as $M\rightarrow\infty$ (to qualify as a delta sequence of 
the Dirichlet type, $M\geq 2$ is required in Eq. (\ref{Lag-dsc})).

{\it Example 4. Interpolative delta sequence kernel}

Let $\{\delta_n\}$ be a sequence of functions 
converging to the delta distribution and 
let $\{x_i \}_0^n$ be $n+1$ zeroes of a
Jacobi polynomial in $(a, b)$.   
\begin{equation}
\Delta_{n}(x,y)=  {\prod_{i=0}^n(x-x_i)  \over (x-y)\prod_{i=0}^n(y-x_i)}
         \sum_{i=0}^n\delta_n(y-x_i), ~~~     x,y \in (a,b)
\end{equation} 
is a delta sequence kernel 
as $n\rightarrow\infty$. This follows from the fact that 
$\int \Delta_{n}(x,y)f(y)dy $ are approximations to the Lagrange 
interpolation formula.

{\it Example 5. The  de la Vall\'{e}e Poussin delta sequence kernel}

The  de la Vall\'{e}e Poussin kernel 
is given by 
\begin{eqnarray}
P_{n,p}(x)&=& {1\over p+1}
\sum_{k=n-p}^nD_k(x)\nonumber\\
&=&
{1\over2\pi}+{1\over\pi}\sum_{k=1}^{n-p}\cos kx +
{1\over\pi}\sum_{k=1}^p
\left[1-{k\over p+1}
\right]\cos[(n-p+k)x] \nonumber\\
&=&
{\sin[(2n+1-p){x\over2}]\sin[(p+1){x\over2}]
 \over 2\pi (p+1) \sin^2({x\over2}) },
   ~~~p=0,\cdots,n;~~n=0,1,\cdots,
\end{eqnarray} 
where $D_k(x)$ are the Dirichlet kernels given by Eq. (\ref{Dirichlet}).
It is interesting to note that  de la Vall\'{e}e Poussin kernel
reduces to the positively defined Fer\'{e}r's kernel $F_{n+1}(x)$  
when $p=n$. 
The de la Vall\'{e}e Poussin delta sequence kernel is given by
\begin{equation}
\delta_{n,p}(x)=\left\{\begin{array}{ll}P_{n,p}(x) 
         &\mbox{for $\mid x\mid\leq \pi $}~~~{\rm for}~~ 
   ~~~p=0,\cdots,n;~~n=0,1,\cdots,\\
       0 &\mbox{otherwise}\end{array}\right.
\end{equation} 
as  $n, p\rightarrow\infty$. The de la Vall\'{e}e Poussin delta
sequence kernel is of Dirichlet type when $p<n$. 

   A simplified de la Vall\'{e}e Poussin kernel given by 
\begin{equation}
\delta_{\alpha}(x)={1\over\pi\alpha} 
{\cos(\alpha x)-\cos(2\alpha x)\over x^2} 
\end{equation} 
is found to be very useful numerically\cite{weijcp99}.

{\it Example 6. DSC kernels 
constructed by orthogonal basis expansions}

Let $\{\psi_i \}$ be a complete orthonormal  $L^2(a,b)$ basis. Then
\begin{equation}\label{orth1}
\delta_{n}(x,y)=  
         \sum_{i=0}^n\psi_i(x)\psi_i(y), ~~~     x,y \in (a,b)
\end{equation} 
are DSC delta sequence kernels. 
In the case of trigonometric functions, we again obtain the 
Dirichlet kernels 
given in the {\it Examples  1} and {\it 3}.  
Hermite function expansion is given by
\begin{equation}\label{HDS}
\delta_{n}(x)=   \exp\left( -{x^2}\right)
       \sum_{k=0}^n \left({-1\over4}\right)^k{1\over \sqrt{\pi} k!}
        H_{2k}(x),
 ~~~  \forall   x\in R
\end{equation} 
where $H_{2k}(x)$
is the usual Hermite polynomial. 
Note that the Hermite's kernel in Eq. (\ref{HDS}) 
has a different form from Eq. (\ref{orth1}). 
This is because it is evaluated at 
$x=0$ in the series expansion.

{\it Example 7. Shannon's  delta sequence}

Shannon's delta kernel (or     
Dirichlet's continuous delta kernel)    
is given by the following (inverse) Fourier    
transform of the characteristic function, $\chi_{[-\alpha,\alpha]}$,    
\begin{eqnarray}\label{Shannon}     
\delta_{\alpha}(x)&=& {1\over 2\pi}      
\int_{-\infty}^{\infty}  \chi_{[-\alpha,\alpha]} e^{-i\xi x}d\xi 
\nonumber\\    
&=& {\sin(\alpha x)\over \pi x}.                                     
\end{eqnarray}    
Alternatively,    
Shannon's delta kernel can be given as an integration    
\begin{equation}    
\delta_{\alpha}(x)=     
{1\over \pi} \int_{0}^{\alpha}\cos(xy) dy,    
\end{equation}     
or as the limit of a continuous product    
\begin{equation}    
\delta_{\alpha}(x)= \lim_{N\rightarrow    
\infty}    
{\alpha\over \pi} \prod_{k=1}^N\cos({\alpha\over2^k}x)    
=\lim_{N\rightarrow    
\infty}    
{1\over 2^N\pi} {\sin(\alpha x)\over \sin({\alpha\over2^N}x)}.    
\end{equation}     
    
Numerically,   
Shannon's delta kernel is the most important, because of   
its property of being an element of the  
Paley-Wiener reproducing kernel Hilbert space $B^2_\pi$.
\begin{equation}\label{eq8}   
f(x)=\int_{-\infty}^{\infty}f(y)   
{\sin\pi(x-y)\over\pi(x-y)}dy   
,~~~~\forall f\in B^2_\pi,     
\end{equation}   
where $\forall f\in B^2_\pi$ indicates that, in its Fourier    
representation, the $L^2$ function $f$ vanishes outside    
the interval $[-\pi,\pi]$.     
The Paley-Wiener reproducing kernel Hilbert space $B^2_\pi$  
is a subspace of the Hilbert space $L^2(R)$.    
It is noted that the reproducing kernel Hilbert space is a special class  
of Hilbert space. For instance,  the space $L^2(R)$ is not    
a reproducing kernel Hilbert space.

{\it Example 8. Generalized Lagrange delta kernel}

Shannon's delta kernel can be derived from the   
generalized Lagrange interpolation formula   
\begin{equation}\label{LD2}  
S_k(x) ={G(x)\over G'(x_k)(x-x_k)},  
\end{equation}  
where $G(x)$ is an entire function given by   
\begin{equation}\label{LD3}  
G(x)=(x-x_0)\prod_{k=1}^{\infty}\left(1-{x\over x_k}\right)  
\left(1-{x\over x_{-k}}\right),  
\end{equation}  
and $G'$ denotes the derivative of $G$.  
For a function bandlimited to $B$, 
the generalized Lagrange  interpolating formula 
$S_k(x)$ of Eq. (\ref{LD2}) can provide an exact result
\begin{equation}\label{eq4-2}  
f(x)=\sum_{k\in Z}f(y_k)S_k(x),  
\end{equation}  
whenever the set of non-uniform sampling points satisfy
\begin{equation}\label{eq4-25}  
\sup_{k\in Z}\left|x_k-{k\pi\over B}\right|< {\pi\over 4B}, 
\end{equation}  
where the symbol $Z$ denotes the set of all integers.   
This is called the Paley and Wiener sampling theorem in the literature.

If $\{ x_k \}_{k\in Z}$ are limited to a set of points on    
a uniform infinite grid ($x_k=k\Delta=-x_{-k}$),  
Eq. (\ref{LD3}) can be simplified   
\begin{eqnarray}\label{LD18}  
G(x)&=&x\prod_{k=-\infty,~k\neq0}^{\infty}  
\left( 1 -{x\over k\Delta}\right)\\  
&=&x\prod_{k=1}^{\infty}  \left(1-{x^2\over k^2\Delta^2}\right)\\  
&=&\Delta{\sin{\pi\over\Delta}x\over \pi}.  
\end{eqnarray}  
Since  $G'(x_k)$ reduces to  
\begin{equation} 
G'(x_k)=(-1)^k  
\end{equation} 
on a uniform grid,   
Eq. (\ref{LD2}) gives rise to  
\begin{eqnarray}\label{LD19}  
S_k(x) &=&{G(x)\over G'(x_k)(x-x_k)}=  
{(-1)^k\sin{\pi\over\Delta}x\over{\pi\over\Delta} (x-k\Delta)}\\  
&=&{\sin{\pi\over\Delta}(x-x_k)\over {\pi\over\Delta}(x-x_k)}\label{LD19-2}.  
\end{eqnarray} 
Obviously,    
${\sin{\pi\over\Delta}(x-x_k)\over {\pi\over\Delta}(x-x_k)}$   
is an approximation to the delta distribution  
\begin{eqnarray}  
\lim_{\Delta \rightarrow 0}   
{\sin{\pi\over\Delta}(x-x_k)\over {\pi\over\Delta}(x-x_k)}  
\rightarrow \delta(x-x_k).  
\end{eqnarray}  
In fact,  
the generalized Lagrange interpolation formula   
directly gives rise to the delta distribution under   
an appropriate limit  
\begin{eqnarray}  
\lim_{\max\Delta x\rightarrow 0}   
S_k(x) =  
\lim_{\max\Delta x\rightarrow 0}   
{G(x)\over G'(x_k)(x-x_k)}  
\rightarrow \delta(x-x_k),  
\end{eqnarray}  
where $\max\Delta x$ is the largest $\Delta x$ on the grid.

\subsection{Connection to wavelets}

The DSC approximation to the delta    
distribution is closely related to the theory of wavelets    
and frames. Mathematically,     
wavelets are functions generated from a single function     
by applying dilation and translation.     
They form building blocks for some spaces,     
such as $L^2(R)$, whether as a frame     
or as an orthonormal basis.    
Such building blocks are computationally important when they     
have certain regularity and localization    
in both time and frequency domains.     
Physically, the wavelet transform is a     
mathematical technique that can be used to split a signal into     
different frequency     
bands or components so that  each component can be     
studied with a resolution matched to its scale, thus providing      
excellent frequency and spatial resolution, and  achieving     
computational efficiency.

Shannon's wavelet is one of the most important examples and its     
scaling function is the Shannon's delta kernel,
\begin{equation}    
\phi(x)= {\sin\pi x \over\pi x}.    
\end{equation}    
As a delta kernel, it is normalized     
\begin{equation}\label{norm}    
\hat{\phi}(0)=\int\phi(x)dx= 1,    
\end{equation}    
and its Fourier transform is given by the characteristic function     
$\hat{\phi}(\omega)= \chi_{[-1/2, 1/2)}$.     
It is easy to see that     
\begin{equation}\label{orth0}    
\sum_{n=-\infty}^{\infty} \hat{\phi}(\omega+n)=1    
\end{equation}    
and     
\begin{equation}\label{orth01}    
\sum_{n=-\infty}^{\infty}\mid \hat{\phi}(\omega+n)\mid^2=1.    
\end{equation}   
Equation (\ref{orth01}) is a consequence of orthonormality.   
In fact,  the sequence of functions 
\mbox{$\{\phi(x-n)\}_{n=-\infty}^{\infty}$} are orthonormal.    
        
Shannon's  mother wavelet can be constructed from the Shannon's   
delta kernel (Shannon's wavelet scaling function)  
\begin{equation}\label{orth1-2}    
\psi(x) ={\sin2\pi x -\sin\pi x \over \pi x},    
\end{equation}    
with its Fourier expression     
\begin{equation}\label{orth2}    
\hat{\psi}(\omega)=\chi_{[-1, 1]}(\omega)-\chi_{[-1/2, 1/2]}(\omega).    
\end{equation}    
This is recognized as the    
ideal band pass filter and it satisfies  the orthonormality conditions    
\begin{equation}\label{orth0-2}    
\sum_{n=-\infty}^{\infty} \hat{\psi}(\omega+n)=1    
\end{equation}    
and     
\begin{equation}\label{orth01-2}    
\sum_{n=-\infty}^{\infty}\mid \hat{\psi}(\omega+n)\mid^2=1.    
\end{equation}     
Technically, it can be shown that a system of     
orthogonal wavelets is generated from a single     
function, a ``mother" wavelet $\psi$, by standard operations of     
translation and dilation    
\begin{equation}\label{FMW}    
\psi_{mn}(x)=    
 2^{-{m\over2}}\psi\left({x\over2^m}-n\right), ~~~m,n\in Z.  
\end{equation}

A family of Shannon's wavelet scaling functions   
$\{\psi_{mn}(x)\}_{n,m\in Z}$   
span a series of orthogonal wavelet subspaces $\{W_m\}_{m\in Z}$ 
satisfying
\begin{equation}   
\bigoplus_{m\in Z}W_m = L^2(R).   
\end{equation}   
Alternatively,  
a family of Shannon's wavelet scaling functions   
$\{\phi_{mn}(x)\}_{n,m\in Z}$ are constructed   
from a single Shannon's delta kernel  
\begin{equation}\label{FMW2}    
\phi_{mn}(x)=    
 2^{-{m\over2}}\phi\left({x\over2^m}-n\right), ~~~m,n\in Z .   
\end{equation}     
They span a series of nested wavelet subspaces 
$\{V_m\}_{m\in Z}$, each corresponds   
to a different resolution   
\begin{equation}  
\cdots \subset V_{-1}\subset V_0\subset V_{1}\subset\cdots\subset L^2(R).  
\end{equation}  
This nested structure provides the conceptual basis for the wavelet   
multiresolution analysis.

From the point of view of signal processing,    
the Shannon's delta kernel   
$\phi_\alpha$ corresponds to a family of     
{\it ideal low pass filters}, each with a different     
bandwidth     
\begin{equation}    
\phi_\alpha (x)= {\sin{\alpha     
x }\over \pi x}.    
\end{equation}    
Their corresponding wavelet expressions,
\begin{equation}\label{orth1-22}    
\psi_{\alpha}(x) ={\sin2\alpha\pi x -\sin\alpha\pi x \over \pi x},    
\end{equation}    
are band pass filters.     
However, the Shannon's wavelet system is seldom used in real  
applications because it requires infinitely many data points.  
In the next subsection, we discuss a practical approach   
for generating powerful filters  from Shannon's delta kernel.

\subsection{Regularization}

Both $\phi(x)$ and its associated wavelet play a crucial    
role in information theory and the theory of   
signal processing. However their    
usefulness is limited by the fact that $\phi(x)$ and $\psi(x)$ are    
infinite impulse response (IIR) filters and their Fourier  
transforms   $\hat{\phi}(\omega)$ and $\hat{\psi}(\omega)$     
are not differentiable. From the computational point of view,  
$\phi(x)$ and $\psi(x)$ do not have     
finite moments in the coordinate space; in other words,     
they are de-localized. This non-local feature in the coordinate   
is related to its bandlimited character in the Fourier   
representation by the Heisenberg uncertainty principle.   

According to the theory of distributions, 
the smoothness, regularity and localization of  
a temper distribution can be improved by a 
function of the Schwartz class.
We apply this principle to regularize singular convolution kernels
\begin{equation}\label{rsinc4-2}    
\tilde{\Phi}_{\sigma}(x) = R_{\sigma}(x)\phi(x),~~~~(\sigma>0),    
\end{equation}     
where $R_{\sigma}$ is a {\it  regularizer} which     
has properties    
\begin{equation}\label{r0}     
\lim_{\sigma\rightarrow\infty}R_{\sigma}(x)=1    
\end{equation}     
and    
\begin{equation}\label{r1}     
R_{\sigma}(0)=1.    
\end{equation}    
Here Eq. (\ref{r0}) is a general condition     
that a  regularizer must satisfy,     
while Eq. (\ref{r1}) is specifically for a {\it delta     
regularizer}, which is used in regularizing a     
delta kernel.    
Various delta regularizers 
 can be used for numerical  
computations. A good example
is the Gaussian 
\begin{equation}\label{r2}     
R_{\sigma}(x)=\exp\left[{-{x^2\over2\sigma^2}}\right].    
\end{equation}  
Gaussian regularizer 
is a Schwartz class function and has
excellent numerical performance.
However, we noted that in certain eigenvalue problems, no 
regularization is required 
if the potential is smooth and bounded from below (e.g., the 
harmonic oscillator potential ${1\over2}x^2$).

Immediate benefit from the regularized Shannon's kernel  
function Eq. (\ref{rsinc4-2}) is that its Fourier transform     
is infinitely differentiable because the Gaussian is  
an element of Schwartz class functions.   
Qualitatively, all kernels of the Dirichlet type 
oscillate in the coordinate representation.  
Shannon's kernel has a long tail which is  
proportional to ${1\over x}$, whereas, 
the regularized kernels decay exponentially fast,  
especially when the $\sigma$ is very small. 
In the Fourier representation, Shannon's kernel  
is an ideal low pass filter, which is discontinuous  
at $\omega={1\over2}$. In contrast, all  
regularized Shannon's kernels have an ``optimal''  
shape in their frequency responses. 
Of course, they all reduce to Shannon's low pass filter at  
the limit 
 \begin{equation}\label{rsinc1-2}    
\lim_{\sigma\rightarrow\infty}\Phi_\sigma(x)=     
\lim_{\sigma\rightarrow\infty}{\sin\pi x\over \pi x}    
e^{-{x^2\over 2\sigma^2}}={\sin\pi x\over \pi x}.  
\end{equation}    
Quantitatively, one can    
examine the normalization of $\Phi_\sigma(x)$       
\begin{eqnarray}\label{rsinc3}    
&&\hat{\Phi}_{\sigma}(0)    
= \int\Phi_{\sigma}( x )dx \nonumber\\    
&&= \sqrt{2\pi}    
{\sigma}    
\sum_{k=0}^{\infty}{(-1)^k\over k!(2k+1)}    
\left({\pi\sigma \over\sqrt{2} }\right)^{2k}.    
\end{eqnarray}    
By means of the error function,     
${\rm erf}(z)={2\over\sqrt{\pi}}\int_{0}^{z}e^{-t^2}dt$,     
Eq. (\ref{rsinc3}) can be rewritten as    
\begin{eqnarray}    
&&\hat{\Phi}_{\sigma}(0)=    
{\rm erf}\left({\pi\sigma\over\sqrt{2}}\right)\nonumber\\    
&&=1-\sqrt{2\over\pi}{1\over\sigma}e^{-{\sigma^2\pi^2\over2 }}    
\int_0^{\infty }e^{-{t^2\over2\sigma^2}-{\pi t}}dt\nonumber\\    
&&=1-{\rm erfc}\left({\pi\sigma\over\sqrt{2}}\right),    
\end{eqnarray}    
where erfc(z) is the complementary error function.    
Note that for a given  $\sigma>0$,      
erfc(${\pi\sigma\over\sqrt{2}}$) is positive definite.    
Thus, $\hat{\Phi}_{\sigma}(0)$     
is always less than unity except at the limit of     
$\sigma\rightarrow\infty$.

In fact, $\Phi_{\sigma}(x)$ does     
not really satisfy the requirement,   
as given by Eq. (\ref{norm}),     
for a wavelet scaling function.  
However, when we choose $\sigma \gg \sqrt{2}/\pi$, which is     
the case in many practical applications, the residue term,     
erfc(${\pi\sigma\over\sqrt{2}}$), approaches     
zero very quickly.   As a result,  
$\hat{\Phi}_{\sigma}(0)$    
is extremely close to unity. Therefore, we call     
the regularized Shannon's delta kernel $\Phi_{\sigma}$     
a quasi-wavelet scaling function.

\subsection{Discretization}

For the purpose of digital computations, it is necessary to  
discretize various delta kernels. To this end, we should  
examine a {\it sampling basis} given by the Shannon's delta 
kernel  
\begin{equation}\label{eq6}   
S_k(x) =K(x,x_k)=    
{\sin\pi(x-x_k)\over\pi(x-x_k)},~~~\forall k\in Z.   
\end{equation}   
This sampling basis is an element of the Paley-Wiener  
reproducing kernel Hilbert space. Hence,  
it provides a {\it discrete representation} of every   
(continuous) function in ${B^2_\pi}$, that is
\begin{equation}\label{eq4}   
f(x)=\sum_{k\in Z}f(x_k)S_k(x),~~~~\forall f\in {B^2_\pi}.   
\end{equation}   
This is recognized as Shannon's sampling theorem and it  
means that one can recover a continuous bandlimited $L^2$  
function from a set of discrete values. 
Equation (\ref{eq4}) is particularly important to   
information theory and the theory of sampling because       
it satisfies the interpolation condition     
\begin{equation}\label{eq5}   
S_n(x_m)=\delta_{n,m}~~,   
\end{equation}   
where $\delta_{m,k}$ is the  Kronecker delta function.     
Note that Shannon's delta kernel is obviously    
interpolative on $Z$. Computationally, being interpolative is    
desirable for numerical accuracy and    
simplicity.    
 
On a grid of arbitrary spacing $\Delta$, Shannon's sampling  
theorem can be modified as 
\begin{equation}\label{eq4-22}   
f(x)=\sum_{k\in Z}f(x_k) 
{\sin{\pi\over\Delta}(x-x_k)\over{\pi\over\Delta}(x-x_k)}    
,~~~~\forall f\in {B^2_{\pi\over\Delta}}.   
\end{equation}   
This suggests that we can  
 discretize the regularized Shannon's delta kernel as       
\begin{equation}\label{rsinc2}    
\Phi_{\sigma,\Delta}(x-x_k)=     
{\sin{\pi\over\Delta}(x-x_k)\over{\pi\over\Delta}(x-x_k)}    
e^{-{(x-x_k)^2\over 2\sigma^2}}.    
\end{equation}    
It is noted that if $\Delta $ is chosen as the spatial     
mesh size (this is, in general,   
not required in signal and image processing),    
$\Phi_{\sigma,\Delta}(x-x_k)$ retains the     
interpolation property,     
\begin{equation}\label{rsinc2-5}    
\Phi_{\sigma,\Delta}(x_m-x_k)= \delta_{m,k}.   
\end{equation}    
This is of particular     
merit for numerical computations.

In practical applications, 
Eq.~(\ref{eq4-22}) can never be realized because it requires  
infinitely many sampling points. Therefore, it is  
both necessary and convenient to    
truncate the infinite summation in Eq. (\ref{eq4-22}) to a finite 
($2M+1$) summation
\begin{equation}\label{eq4-4}   
f(x)\approx\sum_{k=-M}^M   
\delta_{\sigma,\Delta}(x-x_k)f(x_k),   
\end{equation}   
where $\delta_{\sigma,\Delta}(x-x_k)$ is a collective symbol   
for any (regularized) delta kernel. 
The truncation error is dramatically reduced by   
the introduction of a delta regularizer. 
A rigorous proof of this has been 
given by Qian and Wei\cite{qian99}.

The discretization the of the Dirichlet kernel is 
not as straightforward as is for the Shannon's kernel.
However, it can be carried out according to the following 
Dirichlet sampling theorem:

{\bf Theorem} {\em If an $L^2$ function $f(x)$ 
satisfies the Dirichlet boundary 
condition and is periodic in $T$ 
and bandlimited to the highest (radial) frequency 
$2\pi L/T$, it can be exactly reconstructed from a finite
set of $2L+1$ discrete sampling points
\begin{equation}\label{k2-0}  
f(x)=\sum_{k=-L}^{L}f(x_k) 
{\sin\left[{\pi\over\Delta}(x-x_k)\right]\over   
(2L+1)\sin\left[{\pi\over\Delta}{x-x_k\over2L+1}\right]},       
\end{equation} 
where $\Delta=T/(2L+1)$ is the sampling interval
and $x_k=k\Delta$ are the sampling points.}
 
Note that the kernel in Eq. (\ref{k2-0})
differs form that in Eq. (\ref{Dirichlet}). This 
follows from a change in the variable 
$x\rightarrow {\pi\over\Delta y}{2\over2L+1}y$, and   
$\int dx\rightarrow \sum {\pi\over\Delta y}{2\over2L+1}\Delta y$.
The Dirichlet kernel is a {\it reproducing kernel}
for bandlimited $L^2$ periodic functions. Therefore, 
Eq. (\ref{k2-0}) should be the most efficient kernel 
for numerical computations
under the aforementioned conditions. 
However, to facilitate the Dirichlet 
kernel in an unbounded computational domain,  
we use the following regularized discrete expression    
for the Dirichlet kernel:   
\begin{equation}\label{k2}   
{\sin\left[\left(l+{1\over2}\right)(x-x')\right] 
 \over 2\pi \sin[{1\over2}(x-x')]}   
\rightarrow    
{\sin\left[{\pi\over\Delta}(x-x_k)\right]\over   
(2L+1)\sin\left[{\pi\over\Delta}{x-x_k\over2L+1}\right]}   
{\rm exp}\left[{-{(x-x_k)^2\over 2\sigma^2}}\right].   
\end{equation}        
Like the regularized Shannon's kernel filter,   
the present regularized Dirichlet kernel filter has    
rapid decay properties.    
In comparison to  Shannon's kernel, the   
Dirichlet  kernel has one more parameter   
$L$ which can be optimized to achieve better results in    
computations. Usually, we set a sufficiently large $L$ for    
various numerical applications.   
A regularized discrete expression for the   
modified Dirichlet kernel is    
\begin{equation}\label{k3}   
{\sin[\left[l+{1\over2}\right](x-x')] \over 2\pi \tan[{1\over2}(x-x')]}   
\rightarrow    
{\sin\left[{\pi\over\Delta}(x-x_k)\right]\over   
(2L+1)\tan\left[{\pi\over\Delta}{x-x_k\over2L+1}\right]}   
{\rm exp}\left[{-{(x-x_k)^2\over 2\sigma^2}}\right].   
\end{equation}    
Obviously, the regularized Dirichlet kernel
reduces to the regularized Shannon's delta kernel  
when $L$ is sufficiently large  
\begin{eqnarray}\label{k2-2}   
&&\lim_{L\rightarrow \infty}  
{\sin\left[{\pi\over\Delta}(x-x_k)\right]\over   
(2L+1)\sin\left[{\pi\over\Delta}{x-x_k\over2L+1}\right]}   
{\rm exp}\left[{-{(x-x_k)^2\over 2\sigma^2}}\right]  
\nonumber\\  
&=&\lim_{L\rightarrow \infty}  
{\sin\left[{\pi\over\Delta}(x-x_k)\right]\over   
(2L+1)\tan\left[{\pi\over\Delta}{x-x_k\over2L+1}\right]}   
{\rm exp}\left[{-{(x-x_k)^2\over 2\sigma^2}}\right]  
\nonumber\\  
&=&  
{\sin{\pi\over\Delta}(x-x_k)\over{\pi\over\Delta}(x-x_k)}    
e^{-{(x-x_k)^2\over 2\sigma^2}}.    
\end{eqnarray}

The discretization of the de la Vall\'{e}e Poussin kernel   
is given by    
\begin{eqnarray}\label{k4}   
&&{1\over\pi\alpha}{\cos[\alpha (x-x'~)]-\cos[2\alpha (x-x'~)]   
\over (x-x'~)^2}   
\rightarrow \nonumber \\    
&&\hspace*{1cm} {2\over 3} {\cos{\pi\over\bar{\Delta}}(x-x_k)   
-\cos{2\pi\over\bar{\Delta}}(x-x_k)   
\over \left[{\pi\over\bar{\Delta}}(x-x_k)\right]^2 }   
{\rm exp}\left[{-{(x-x_k)^2\over 2\sigma^2}}\right],   
\end{eqnarray}    
where $\bar{\Delta}={3\over2}\Delta$.   
Since $\pi/\Delta$ is proportional to the highest frequency   
which can be reached in the Fourier representation,   
the $\Delta$ should be very small for a given problem involving    
very oscillatory functions or very high frequency components.

It is noted that by definition, 
the Lagrange interpolation   formula   
\begin{equation} \label{L1}  
L_{M,k}(x)=  
\prod_{i=k-M,i\neq k}^{k+M}{x-x_i\over x_k-x_i}  
\end{equation}
is already discretized. However,   
its regularized forms   
\begin{equation} \label{L2}  
\delta_{\sigma}(x-x_k)=  
\left[\prod_{i=k-M,i\neq k}^{k+M}{x-x_i\over x_k-x_i}\right]  
\exp\left(-{(x-x_k)^2\over 2\sigma^2}\right)  
\end{equation}   
are very good low pass filters.  
Unlike the generalized Lagrange interpolation   
formula, both expressions (\ref{L1}) and (\ref{L2}) are  
compactly supported kernels. As low pass filters, they require only   
only a finite number of signals and their Fourier transforms  
have smoothened  shoulders as those of regularized Shannon's  
delta kernels.

\subsection{Approximation of derivatives}

For the solution of differential equations, an  
approximation to  the derivatives  is required.  
Such an approximation can be constructed by using 
DSC kernels of the delta type with $n\neq0$.  
Let us consider a one-dimensional, $n$th order     
DSC kernel of the delta type   
\begin{equation}
\delta^{(n)}_{\sigma,\Delta}(x-x_k),~~~~(n=0,1,2,\cdots).
\end{equation}
Here $\delta^{(0)}_{\sigma,\Delta}(x-x_k)   
=\delta_{\sigma,\Delta}(x-x_k)$ is the DSC delta kernel 
described in    Eq. (\ref{eq4-4}).     
The higher order derivative terms   
$\delta^{(n)}_{\sigma,\Delta}(x_m-x_k)$ are given by    
differentiation   
\begin{equation}\label{deriv}   
\delta^{(n)}_{\sigma,\Delta}(x_m-x_k)=\left[   
\left({d\over dx}\right)^n \delta_{\sigma,\Delta}(x-x_k)\right]_{x=x_m}.   
\end{equation} 
These derivatives can be regarded as high pass filters. 
The filters corresponding to the  
derivatives of Shannon's kernel decay slowly as  
$x$ increases, whereas, the regularized filters are Schwartz class  
functions and have controlled 
residual amplitudes at large $x$ values.   
In the Fourier representation, the 
derivatives of Shannon's kernel are  
discontinuous at certain points. In contrast,  
all the derivatives of regularized kernels are  
continuous and can be made very close to those of Shannon's, 
if desired.

The differentiation in Eq. (\ref{deriv}) can  
be {\it analytically}              
carried out for a given $\delta_{\sigma,\Delta}(x-x_k)$.   
For example, if $\delta_{\sigma,{\Delta}}(x-x_k)=   
{\sin{\pi\over\Delta}(x-x_k)\over{ \pi\over\Delta }(x-x_k)}   
e^{-{(x-x_k)^2\over 2\sigma^2}}$, we have for $x\neq x_k$    
\begin{eqnarray}\label{deriv1}   
\delta^{(1)}_{\sigma,{\Delta}}(x-x_k)&=&   
{\cos{\pi\over\Delta}(x-x_k)   
\over (x-x_k)}\exp\left(-{(x-x_k)^2\over2\sigma^2}\right)   
\nonumber \\ &-&   
{\sin{\pi\over\Delta}(x-x_k)\over{\pi\over\Delta} (x-x_k)^2}    
\exp\left(-{(x-x_k)^2\over2\sigma^2}\right)   
\nonumber \\ &-&   
{\sin{\pi\over\Delta}(x-x_k)\over   
{\pi\over\Delta}\sigma^2}\exp\left(-{(x-x_k)^2\over2\sigma^2}\right), 
\end{eqnarray}   
\begin{eqnarray}\label{deriv2}   
\delta^{(2)}_{\sigma,{\Delta}}(x-x_k)   
&=&   
-{{\pi\over\Delta}\sin{\pi\over\Delta}(x-x_k)   
\over (x-x_k)}    
\exp\left(-{(x-x_k)^2\over2\sigma^2}\right)   
\nonumber \\    
&-&2{\cos{\pi\over\Delta}(x-x_k)   
\over (x-x_k)^2}   
\exp\left(-{(x-x_k)^2\over2\sigma^2}\right)   
\nonumber \\    
&-&2{\cos{\pi\over\Delta}(x-x_k)   
\over \sigma^2}   
\exp\left(-{(x-x_k)^2\over2\sigma^2}\right)   
\nonumber \\ &+&   
2{\sin{\pi\over\Delta}(x-x_k)\over   
{\pi\over\Delta}(x-x_k)^3 }   
\exp\left(-{(x-x_k)^2\over2\sigma^2}\right)   
\nonumber \\ &+&   
{\sin{\pi\over\Delta}(x-x_k)\over   
{\pi\over\Delta}(x-x_k)\sigma^2 }   
\exp\left(-{(x-x_k)^2\over2\sigma^2}\right)   
\nonumber \\ &+&   
{\sin{\pi\over\Delta}(x-x_k)\over   
{\pi\over\Delta} \sigma^4 }(x-x_k)   
\exp\left(-{(x-x_k)^2\over2\sigma^2}\right),   
\end{eqnarray}   
   
\begin{eqnarray}\label{deriv3}   
\delta^{(3)}_{\sigma,{\Delta}}(x-x_k)   
&=&   
-{{\pi^2\over\Delta^2}\cos{\pi\over\Delta}(x-x_k)   
\over (x-x_k)}    
\exp\left(-{(x-x_k)^2\over2\sigma^2}\right)   
\nonumber \\    
&+&3{{\pi\over\Delta}\sin{\pi\over\Delta}(x-x_k)   
\over (x-x_k)^2}   
\exp\left(-{(x-x_k)^2\over2\sigma^2}\right)   
\nonumber \\   
&+&3{{\pi\over\Delta}\sin{\pi\over\Delta}(x-x_k)   
\over \sigma^2}   
\exp\left(-{(x-x_k)^2\over2\sigma^2}\right)   
\nonumber \\    
&+&6{\cos{\pi\over\Delta}(x-x_k)   
\over (x-x_k)^3 }   
\exp\left(-{(x-x_k)^2\over2\sigma^2}\right)   
\nonumber \\    
&+&   
3{\cos{\pi\over\Delta}(x-x_k)\over   
(x-x_k) \sigma^2 }   
\exp\left(-{(x-x_k)^2\over2\sigma^2}\right)   
\nonumber \\    
&+&   
3{(x-x_k)\cos{\pi\over\Delta}(x-x_k)\over   
\sigma^4 }   
\exp\left(-{(x-x_k)^2\over2\sigma^2}\right)   
\nonumber \\    
&-&   
6{\sin{\pi\over\Delta}(x-x_k)\over   
{\pi\over\Delta}(x-x_k)^4}   
\exp\left(-{(x-x_k)^2\over2\sigma^2}\right)   
\nonumber \\    
&-&   
3{\sin{\pi\over\Delta}(x-x_k)\over   
{\pi\over\Delta}(x-x_k)^2\sigma^2}   
\exp\left(-{(x-x_k)^2\over2\sigma^2}\right)   
\nonumber \\    
&-&   
{(x-x_k)^2\sin{\pi\over\Delta}(x-x_k)\over   
{\pi\over\Delta}\sigma^6}   
\exp\left(-{(x-x_k)^2\over2\sigma^2}\right),   
\end{eqnarray}   
and    
\begin{eqnarray}\label{deriv4}   
\delta^{(4)}_{\sigma,{\Delta}}(x-x_k)   
&=&   
4{{\pi^2\over\Delta^2}\cos{\pi\over\Delta}(x-x_k)   
\over (x-x_k)^2}    
\exp\left(-{(x-x_k)^2\over2\sigma^2}\right)   
\nonumber \\    
&+&   
{{\pi^3\over\Delta^3}\sin{\pi\over\Delta}(x-x_k)   
\over (x-x_k)}    
\exp\left(-{(x-x_k)^2\over2\sigma^2}\right)   
\nonumber \\    
&+&   
4{{\pi^2\over\Delta^2}\cos{\pi\over\Delta}(x-x_k)   
\over \sigma^2 }   
\exp\left(-{(x-x_k)^2\over2\sigma^2}\right)   
\nonumber \\   
&-&12{{\pi\over\Delta}\sin{\pi\over\Delta}(x-x_k)   
\over  (x-x_k)^3 }   
\exp\left(-{(x-x_k)^2\over2\sigma^2}\right)   
\nonumber \\    
&-&6{{\pi\over\Delta} \sin{\pi\over\Delta}(x-x_k)   
\over (x-x_k)\sigma^2 }   
\exp\left(-{(x-x_k)^2\over2\sigma^2}\right)   
\nonumber \\    
&-&6{{\pi\over\Delta}(x-x_k)\sin{\pi\over\Delta}(x-x_k)   
\over \sigma^4 }   
\exp\left(-{(x-x_k)^2\over2\sigma^2}\right)   
\nonumber \\    
&-&   
24{\cos{\pi\over\Delta}(x-x_k)\over   
(x-x_k)^4 }   
\exp\left(-{(x-x_k)^2\over2\sigma^2}\right)   
\nonumber \\    
&-&   
12{\cos{\pi\over\Delta}(x-x_k)\over   
(x-x_k)^2 \sigma^2}   
\exp\left(-{(x-x_k)^2\over2\sigma^2}\right)   
\nonumber \\    
&-&   
4{(x-x_k)^2\cos{\pi\over\Delta}(x-x_k)\over   
\sigma^6 }   
\exp\left(-{(x-x_k)^2\over2\sigma^2}\right)   
\nonumber \\    
&+&   
24{\sin{\pi\over\Delta}(x-x_k)\over   
{\pi\over\Delta}(x-x_k)^5}   
\exp\left(-{(x-x_k)^2\over2\sigma^2}\right)   
\nonumber \\    
&+&   
12{\sin{\pi\over\Delta}(x-x_k)\over   
{\pi\over\Delta}(x-x_k)^3\sigma^2}   
\exp\left(-{(x-x_k)^2\over2\sigma^2}\right)   
\nonumber \\    
&+&   
3{\sin{\pi\over\Delta}(x-x_k)\over   
{\pi\over\Delta} (x-x_k)\sigma^4}   
\exp\left(-{(x-x_k)^2\over2\sigma^2}\right)   
\nonumber \\    
&-&   
2{(x-x_k)\sin{\pi\over\Delta}(x-x_k)\over   
{\pi\over\Delta}\sigma^6}   
\exp\left(-{(x-x_k)^2\over2\sigma^2}\right)   
\nonumber \\    
&+&   
{(x-x_k)^3\sin{\pi\over\Delta}(x-x_k)\over   
{\pi\over\Delta}\sigma^8}   
\exp\left(-{(x-x_k)^2\over2\sigma^2}\right)   
\end{eqnarray}   
At $x=x_k$, it is convenient to evaluate these derivatives   
separately   
\begin{eqnarray}\label{limit}   
\delta^{(1)}_{\sigma,{\Delta}}(0)&=&0   
\\   
\delta^{(2)}_{\sigma,{\Delta}}(0)&=&   
-{1\over3}{3+ {\pi^2\over\Delta^2}\sigma^2\over \sigma^2}   
\\    
\delta^{(3)}_{\sigma,{\Delta}}(0)&=&0   
\end{eqnarray}   
and   
\begin{eqnarray}\label{limit2}   
\delta^{(4)}_{\sigma,{\Delta}}(0)=   
{1\over5}{15+ 10{\pi^2\over\Delta^2}\sigma^2   
+{\pi^4\over\Delta^4}\sigma^4\over \sigma^4}.
\end{eqnarray}   
Similar expressions for other DSC kernels described in the 
last subsection can be easily derived. 
The performance of a few DSC kernels for fluid dynamic 
computations and structural analysis was given in 
Ref. \cite{weicmame20}.
Note that the differentiation matrix in Eq. (\ref{deriv})   
is generally banded. This has a distinct
advantage in large scale computations.   
   
For numerical computations, it turns out that the approximate
reproducing kernel has much less truncation errors
for interpolation and numerical differentiation.
Qian and the present author\cite{qian99}  have recently 
given the following theorem for truncation errors.\\
{\bf Theorem} {\it Let $f$ be a function $f\in L^2(R)\cap C^{s}(R)$ 
and bandlimited to 
$B$, $(B< \frac{\pi}{\Delta},~ \Delta$ is the grid spacing).  
For a fixed $t\in R $ 
and $\sigma>0$, denote $g(x)=f(x)H_k(\frac{t-x}{\sqrt{2}\sigma})$, 
where $H_k(x)$ is the
$k$th order Hermite polynomial. If $g(x)$ satisfies 
\begin{equation} 
g'(x)\leq
g(x)\frac{(x-t)}{{\sigma}^2} 
\end{equation} 
for $x\geq t+(M_1-1)\Delta$, and
\begin{equation} 
g'(x)\geq g(x)\frac{(x-t)}{{\sigma}^2} 
\end{equation} 
for $x\leq
t-M_2\Delta$, where $M_1,M_2\in \mathcal{N}$,  then for any 
$s\in \mathcal{Z}^{+}$
\begin{eqnarray} 
&&\left\|
f^{(s)}(t)-\sum_{n=\lceil\frac{t}{\Delta}\rceil-M_2}^{\lceil
\frac{t}{\Delta}
\rceil+M_1}f(n\Delta)
\left[\frac{\sin\frac{\pi}{\Delta}(t-n\Delta)}{\frac{\pi}{\Delta}
(t-n\Delta)}
\exp(-\frac{(t-n\Delta)^2}{2{\sigma}^2})\right]^{(s)}
\right\|_{L^2(R)}\nonumber\\ 
&&\leq
\sqrt{3}\left[\frac{\| f^{(s)}(t)\|_{L^2(R)}}{2\pi\sigma
(\frac{\pi}{\Delta}-B)
\exp(\frac{{\sigma}^2(\frac{\pi}{\Delta}-B)^2}{2})}\right.
\nonumber\\ 
&&\left.+ \frac{\|
f(t)\|_{L^2(R)}\sum_{i+j+k=s}\frac{s!{\pi}^{i-1}
H_k(\frac{-M_1\Delta}{\sqrt{2}\sigma})}
{i!k!{\Delta}^{i-1}(\sqrt{2}\sigma)^k((M_1-1)\Delta)^{j+1}}}
{\exp(\frac{(M_1\Delta)^2}{2{\sigma}^2})} \right.\nonumber\\ 
&&\left.+\frac{\|
f(t)\|_{L^2(R)}\sum_{i+j+k=s}\frac{s!{\pi}^{i-1}
H_k(\frac{-M_2\Delta}{\sqrt{2}\sigma})}
{i!k!{\Delta}^{i-1}(\sqrt{2}\sigma)^k(M_2\Delta)^{j+1}}}
{\exp(\frac{(M_2\Delta)^2}{2{\sigma}^2})}\right],
\label{eq30} 
\end{eqnarray}
where superscript, $(s)$, denotes the $s$th order derivative.}

The proof and detailed discussion 
(including a comparison with the truncation errors of 
Shannon's sampling theorem) are given in Ref. \cite{qian99} 
and are beyond the scope of this paper.
This theorem provides a guide to 
the choice of $M$, $\sigma$ and $\Delta$.
For example, in the case of interpolation ($s=0$), if the 
$L_2$ norm error is set to $10^{-\eta}$ ($\eta>0$), 
the following relations can be deduced from Eq. (\ref{eq30})
\begin{eqnarray}
r(\pi-B\Delta)>\sqrt{4.61\eta},
\end{eqnarray}
and
\begin{eqnarray}
{M\over r}>\sqrt{4.61\eta},
\end{eqnarray}
where $r=\sigma/\Delta$ (The choice of $\sigma$ is always 
proportional to 
$\Delta$ so that the width of the Gaussian envelope
varies with the central frequency). 
The first inequality states that for a given grid 
size $\Delta$, a large $r$ is required 
for approximating high frequency 
component of an $L^2$ function.
The second inequality indicates that 
if one chooses the ratio $r=3$, then the 
half bandwidth $M\sim30$ can be used to ensure the highest 
accuracy in a double precision computation ($\eta=15$). 
However, for lower accuracy requirement, a much 
smaller half bandwidth can be used.
In general, the value of $r$ is proportional to $M$.
An appropriate value of $M$ is determined by the accuracy
requirement. 
This theoretical estimation is in excellent agreement with 
an  earlier numerical test\cite{weicpl}.

\subsection{Treatment of boundary conditions}

While solving a set  of differential equations, 
boundary conditions should be satisfied.
In a global method, the kernel must be 
constructed in an adaptive manner near the boundary. 
The major drawback of such an algorithm is its difficulty in
the construction of adaptive kernels for problems involving
complex geometries and boundary conditions. Thus, global methods     
have been relatively less successful in dealing with these 
problems compared to their advantage
 in solving problems with  simple geometries
and boundary conditions. In contrast, the DSC algorithm utilizes
a completely different philosophy for the kernel construction
--- the differentiation kernel is 
the same everywhere and  is translation
invariant on the grid. Therefore, the DSC algorithm is
very flexible in dealing with a variety of boundary conditions 
and geometries. 
Moreover, the DSC treatment of boundary condition has its 
mathematical justification. In fact,  for
a continuous function,  a
derivative at a point exists if and only if 
both the left and the right derivatives 
exist and are equal. Therefore, with a finite 
computational domain, a boundary 
condition involving differentiation (such as 
the Neumann boundary condition) does not make 
sense. This is an ill-posed problem from the point of view
of mathematical differentiation. The DSC approach to this 
problem is to
extend the domain of definition for the system so 
that the ``boundary condition'' is well defined (i.e., 
differentiation is well-defined right on the boundary).   
The philosophy behind DSC is that, the original 
singular convolution has 
to be recovered from the DSC at the 
limit of $\Delta x\rightarrow 0$, everywhere 
in the computational domain, including boundaries. 
Therefore, at a boundary, a fictitious domain is required 
to ensure that the boundary condition is exactly satisfied 
at the continuous limit.

In an explicit treatment, boundary conditions are 
easily implemented in the DSC algorithm
by appropriate boundary extensions, which were discussed
in Ref.\cite{weijcp99}.
In the present work, we consider implicit cases, where the 
boundary conditions are to be satisfied by a
set of linear algebraic equations.

We assume that a general boundary condition at the left boundary
is given by
\begin{equation}\label{boundary}
\sum_{n=0}^N K_nf^{(n)}(x_0)=0,
\end{equation}
where $K_n$ is a constant and $x_0$ is the boundary point.
In the DSC treatment, we make an assumption for the
relation between the 
inner nodes and the outer nodes on the left boundary
\begin{equation}
f(x_{-i})-f(x_0)=a_i\left[ (f(x_i)-f(x_0)\right],  
\end{equation}
where parameters $a_i, i=1,\ldots, M$ are to be 
determined.    After a rearrangement, we obtain 
\begin{equation}
f(x_{-i})=a_if(x_i) +(1-a_i)f(x_0), ~i=1, 2, \ldots, M. 
\end{equation}
  
According to Eq. (\ref{deriv}), we approximate
the $n$th derivative and second derivative at the left 
boundary by
\begin{eqnarray} 
f^{(n)}(x_0)&=&\sum_{i=-M}^MC^n_if(x_i)\nonumber\\ 
&=&C^n_0 f(x_0)+\sum_{i=1}^M C^n_i\left[ 
\left(1+(-1)^na_i\right)f(x_i)+(-1)^n(1-a_i)f(x_0)\right],
\end{eqnarray}
where $C^n_i$ are
coefficients given by a DSC kernel, e.g.,  
$C^n_i=\delta^{(n)}_{\sigma,\Delta}(x_0-x_i)$. 
Therefore, the boundary condition (\ref{boundary}) on the left 
boundary is given by 
\begin{eqnarray}\label{boundary2}
&&K_0f(x_0)+
\sum_{n=1}^NK_n\left[ C_0^nf(x_0) \right. \nonumber\\ 
&&\left. +\sum_{i=1}^MC_i^n(-1)^n(1-a_i)f(x_0)\right.\nonumber\\ 
&&\left.  +\sum_{i=1}^MC_i^n(1+(-1)^na_i)f(x_i)
\right]=0.
\end{eqnarray}
In general, Eq. (\ref{boundary2}) may not have an exact solution for implicit 
schemes. However, 
for a class of boundary conditions
\begin{equation}
f(x_0)=0, ~~~~ 
\sum_{n=1}^N K_nf^{(n)}(x_0)=0,
\end{equation}
the exact solution is given by 
\begin{equation}\label{boundary3}
a_i=\frac{\sum_{l=1}^{N\over2} 
K_{2l}C^{2l}_i-K_{2l-1}C_i^{2l-1}}
{\sum_{l=1}^{N\over2} 
K_{2l-1}C^{2l-1}_i-K_{2l}C_i^{2l}}
~~~i=1, 2, \ldots, M. 
\end{equation}
The following a few special cases are important 
in practical applications.

One case is the clamped edges, i.e., the boundary 
conditions require 
\begin{equation}
f(x_0)=0, ~~~~ f^{(1)}(x_0)=0. 
\end{equation}
These are satisfied by choosing $a_i=1, i=1,2,\ldots,M$.  
This is the so called {\it symmetric extension}\cite{weijcp99}.

Another case is simply supported edges, i.e. the boundary 
conditions (\ref{boundary})  reduce to
\begin{equation} 
f(x_0)=0,~~~~ f^{(2)}(x_0)=0.
\end{equation}
These are satisfied by choosing $a_i=-1, i=1,2,\ldots,M$. 
This is the so called 
{\it anti-symmetric extension}\cite{weijcp99}.

The other special case of Eq. (\ref{boundary})  
is called transversely supported edge 
and the boundary conditions are given by 
\begin{equation}
f(x_0)=0, ~~~~ f^{(2)}(x_0)+K_1 f^{(1)}(x_0)=0.
\end{equation}
In this case, Eq. (\ref{boundary2}) is reduced to,
\begin{equation}\label{b2}
\sum_{i=1}^M [(1+a_i) C_i^2+K_1 (1-a_i) C_i^1] f(x_i)=0. 
\end{equation}
One way to satisfy Eq. (\ref{b2}) is to choose
\begin{equation}\label{boundary4}
a_i=\frac{K_1C_i^1+ C_i^2}{K_1C_i^1- C_i^2}, ~i=1, 2, \ldots, M. 
\end{equation}
Expressions for the right boundary can be derived in 
a similar way.

\section{Unified Features For Solving Differential Equations}

To solve a differential equation, 
one can start, either by 
approximating the  original differential 
operator or by approximating 
the actual solution of the differential 
equation  while maintaining 
the original differential operator. 
The latter is accomplished by 
explicitly defining a functional form for 
approximations.
Let us assume that the differential 
equation has the general form
\begin{equation}\label{Eq1}
  {\cal  L} u({ x})  = f({ x}),   ~~ { x}\in \Omega,
\end{equation}
where ${\cal  L}$ is a linear operator and 
$u({ x})$ is the unknown solution of interest. Here 
$f(x)$ is a known source term, 
$\Omega$ denotes the domain over which the 
differential equation applies.

       The approximate solution is sought from a finite 
set of $N$ DSC trial functions of a given 
{\it resolution} $\alpha$, 
denoted by $S^{N,M}_{\alpha,\sigma}$ with $M$ being 
the half width of support of each element. 
Here $\sigma$ is a {\it regularization} parameter
for improving the regularity of the set.  
The case of regularization free is easily obtained by setting 
$\sigma\rightarrow\infty$.
Elements of the set $S^{N,M}_{\alpha,\sigma}$ can be explicitly 
given by $\{\phi^M_{\alpha,\sigma;1}, \phi^M_{\alpha,\sigma;2},...,
\phi^M_{\alpha,\sigma; N}\}$. For a given computational 
domain, the resolution parameter $\alpha$ is determined 
by $N$.

    An important property of the DSC trial functions 
$\{\phi^M_{\alpha,\sigma;k}\}$
is that when the  trial function is free of regularization,
each  member of the set is a {\it reproducing kernel} 
at highest resolution
\begin{equation}\label{prop1}
\lim_{\alpha\rightarrow\infty}
<\phi^M_{\alpha,\sigma; k},\eta >=\eta({ x}_k),
\end{equation}
where $< \cdot,\cdot >$ denotes the standard inner product.
In fact, if an appropriate 
basis is used for $\phi$ and the limit on $\sigma$ is
taken, $\phi$ of each resolution 
can be a reproducing kernel for $L^2$ functions 
bandlimited to appropriate sense. In general, we require
the low pass filter property 
that for given $\alpha\neq0,\sigma\neq0$ and $M\gg 0$ 
\begin{equation}\label{prop2}
<\phi^M_{\alpha,\sigma; k},\eta >\approx\eta({ x}_k).
\end{equation}
This  converges uniformly
when the resolution is refined, e.g.,
$\alpha\rightarrow\infty$. 
A few examples of such DSC trial functions are given in 
Refs. \cite{weijcp99} and \cite{weicpl}, and many more  examples
are constructed in the previuos section.
Equations (\ref{prop1}) and (\ref{prop2}) are special 
requirements satisfied by the DSC kernels of delta 
type\cite{weijcp99}.

    In the present DSC approach, an approximation 
to the function of interest $u(x)$ can be expressed 
as a linear combination   
\begin{equation}\label{Eq2-3}
U^{N,M}_{\alpha,\sigma}({ x}) =\sum_{k=1}^{N}
U_{\alpha,\sigma; k} 
\phi^M_{\alpha,\sigma; k}({ x}),
\end{equation}
where ${x}$ is an independent variable and
$U_{\alpha,\sigma; k}$ 
is the desired DSC approximation to the solution at point $x_k$.
This structure is due to 
the DSC trial function property (\ref{prop2}) and it
dramatically 
simplifies the solution procedure in
practical computations.

In this formulation, we choose the set 
$S^{N,M}_{\alpha,\sigma}$ a priori, 
and then determine the 
coefficients $\{ U_{\alpha,\sigma; k} 
\}$ 
so that $U^{N,M}_{\alpha,\sigma}({ x})$ is a good 
approximation to $u({ x})$. 
To determine $U_{\alpha,\sigma; k}$, we minimize the 
amount by which $U^{N,M}_{\alpha,\sigma}({ x})$ 
fails to satisfy the original governing 
equation (\ref{Eq1}). A measure of 
this discrepancy can be  defined as
\begin{equation}\label{Eq3}
R^{N,M}_{\alpha,\sigma}({ x}) \equiv  
{\cal L} U^{N,M}_{\alpha,\sigma}({ x}) - f({ x}),
\end{equation}
where $R^{N,M}_{\alpha,\sigma}({ x})$ is the residual 
for a particular choice of resolution, regularization and
half width of the support. Note that Eq. (\ref{Eq3})
is similar to the usual statement in the method of 
weighted residuals. However,  the approximation 
$U^{N,M}_{\alpha,\sigma}({ x})$ is 
constructed by using the DSC trial functions,
$\phi^M_{\alpha,\sigma; k}({ x})$, in the present 
treatment. Let Eq. (\ref{Eq1}) and its associated 
boundary conditions be well-posed, 
then there exists 
a unique solution $u({ x})$ which generally resides 
in an infinite-dimensional space. 
Since the DSC approximation $U^{N,M}_{\alpha,\sigma}$ is 
constructed from a finite-dimensional set, 
it is generally the case 
that $U^{N,M}_{\alpha,\sigma}({ x})\neq u({ x})$ and 
therefore $R^{N,M}_{\alpha,\sigma}({ x}) \neq 0$.

{\it Galerkin.}~
We seek to optimize $R^{N,M}_{\alpha,\sigma}({ x})$
 by forcing it to zero in 
a {\it weighted average sense} over the domain 
$\Omega$. A convenient starting point is 
the Galerkin 
\begin{equation}\label{Eq4}
\int_\Omega R^{N,M}_{\alpha,\sigma}({ x}) 
\phi^{M'}_{\alpha',\sigma'; l} ({ x}) d{ x} = 0,~~
\phi^{M'}_{\alpha',\sigma'; l}({ x}) 
\in S^{N',M'}_{\alpha',\sigma'}, 
\end{equation}
where the  weight set 
$S^{N',M'}_{\alpha',\sigma'}$  can be simply 
chosen being identical to the DSC trial function set
$S^{N,M}_{\alpha,\sigma}$. We refer to
Eq. (\ref{Eq4}) as a DSC-Galerkin statement.

{\it Collocation.}~
First, we note that in view of Eq. (\ref{prop1}), 
the present DSC-Galerkin statement reduces to 
a collocation one at the limit of $\alpha'$
\begin{equation}\label{Eq5}
\lim_{\alpha'\rightarrow\infty}
\int_\Omega R^{N,M}_{\alpha,\sigma}({ x}) 
\phi^{M'}_{\alpha',\sigma'; l} ({ x}) d{ x}=
R^{N,M}_{\alpha,\sigma}({ x}_l)=0,
\end{equation}     
where $\{{  x}_l \}$ is the set of collocation points.
However, in digital computations, we cannot take the above
limits. It follows from the low pass filter property of the 
DSC trial functions, Eq. (\ref{prop2}), that 
\begin{equation}\label{Eq6}
\int_\Omega R^{N,M}_{\alpha,\sigma}({ x}) 
\phi^{M'}_{\alpha',\sigma'; l} ({ x}) d{ x}\approx
R^{N,M}_{\alpha,\sigma}({ x}_l)\approx0.
\end{equation}     
It can be proven that for an appropriate choice of 
$S^{N',M'}_{\alpha',\sigma'}$,
the first approximation of Eq. (\ref{Eq6}) 
converges uniformly. 
The difference between
the true DSC-collocation, 
\begin{equation}\label{Eq5-2}
\lim_{\alpha'\rightarrow\infty}
R^{N,M}_{\alpha,\sigma}({ x}_l)=0,
\end{equation}     
and the {\it Galerkin induced collocation}, (\ref{Eq6}), 
diminishes to zero for 
appropriate DSC trial functions.

{\it Global and local.}~
Global approximations to a function and its derivatives are realized
typically by a set of truncated $L^2(a,b)$ function expansions.
It is called global because the values of a function and its
derivatives at a particular point $x_i$ in the coordinate
space involve the {\it full} set of grid points in a computational 
domain $\Omega$. Whereas a local method does so by requiring only 
a few neighborhood points.
In the present DSC approach, since the choices of $M$ and/or $M'$ are 
independent of $N$, one can choose $M$ and/or $M'$
so that a function and its
derivatives at a particular point $x_l$ are approximated either
by the full set of grid points in the computational 
domain $\Omega$ or just by a few grid points in the 
neighborhood.
In fact, this freedom for the selection of $M$ endows the DSC 
algorithm  with {\it controllable accuracy} for solving 
differential equations and the flexibility in handling 
complex geometries.

{\it Finite Difference.}~
In the finite difference method, the differential 
operator is approximated by difference operations.
In the present approach, the DSC-collocation 
expression of Eq. (\ref{Eq6})
is equivalent to a $2M+1$ (or $2M$) term  finite difference method.
This follows from the fact that the DSC approximation to the 
$n$th order derivative of a function can be rewritten as 
\begin{equation}\label{Eq10}
\left.{d^q u\over dx^q}\right|_{x=x_k}
\approx
\sum_{l=k-M}^{k+M}c^{q}_{kl,M}u(x_l),
\end{equation}     
where $c^{q}_{kl,M}$ are a set of DSC weights for the finite
difference approximation and are given by 
\begin{equation}\label{Eq11}
c^{q}_{kl,M}= 
\left.{d^q \over dx^q}\phi^{M}_{\alpha,\sigma; l} ({ x})
\right|_{{ x}={ x}_k}. 
\end{equation}     
Note that the standard finite difference scheme is 
obtained if $\phi^{M}_{\alpha,\sigma; l}$
is chosen as the Lagrange kernel as discussed earlier.
Obviously, for each different choice of 
$\phi^M_{\alpha,\sigma}$, we have a different 
DSC-finite difference approximation. 
Hence, the present DSC approach is a generalized finite 
difference method. This DSC-finite difference was tested 
in previous studies\cite{weicpl}.
When $M=1$, the DSC-finite difference
approximation reaches its lowest order 
limit and the resulting matrix is tridiagonal.
In this case,
by appropriately choosing the parameter $\sigma$,
the present DSC weights $c^{q}_{kl,M}$
can always be made exactly the same as those of  
the standard second order central difference scheme, 
i.e. ${1\over2\Delta}, 0, -{1\over2\Delta}$ for the first order 
derivative and ${1\over\Delta^2}, -{2\over\Delta^2}, {1\over\Delta^2}$
for the second order derivative. Here $\Delta$ is the grid spacing.
However,  for a given numerical bandwidth,
the DSC-finite difference approximation
does not have to be the same as the standard 
finite difference scheme 
and can be optimized in a practical 
application by varying $\sigma$. 
Another important choice of the DSC bandwidth is that 
$M=N$, where $N$ is the matrix length. Obviously, the 
computational matrix is no longer banded and this is a case 
we called a ``global finite difference method''.  
These interesting features are illustrated by using 
a few numerical examples in the next section.

\section{Application to Electromagnetics}          

In this section, we examine the usefulness, test 
the accuracy and explore the limitation of the DSC
approach for solving problems in electromagnetics. 
Many DSC kernels discussed in the previous sections
are suitable for being  used as DSC trial functions.
For simplicity, we focus on three DSC kernels of 
the Dirichlet type, a regularized Shannon's kernel (RSK), 
\begin{equation}\label{K1}
\phi^M_{{\pi\over\Delta},\sigma;k}(x)=
{\sin{\pi\over\Delta}(x-x_k)\over
{\pi\over\Delta}(x-x_k)}
\exp\left[ {-{(x-x_k)^2\over2\sigma^2}}\right],
\end{equation}     
a regularized Dirichlet kernel (RDK), 
\begin{equation}\label{K2}
\phi^M_{{\pi\over\Delta},\sigma;k}(x)=
{\sin\left[{\pi\over\Delta}(x-x_k)\right]\over
(2L+1)\sin\left[{\pi\over\Delta}{x-x_k\over2L+1}\right]}
\exp\left[ {-{(x-x_k)^2\over2\sigma^2}}\right],
\end{equation} 
and a regularized Lagrange kernel (RLK)
\begin{equation}\label{K3}
\phi^M_{{\pi\over\Delta},\sigma;k}(x)=
\prod_{i=k-L,i\neq k}^{i=k+L}{x-x_i\over x_k-x_i}
\exp\left[ {-{(x-x_k)^2\over2\sigma^2}}\right],
\end{equation} 
for our numerical experiments.
Note that the resolution is given by 
$\alpha={\pi\over\Delta}$, which is the frequency bound in 
the Fourier representation. 
The half bandwidth, $M$,  
can be chosen to interplay between the local limit and the global
limit and is set to 40 in all calculations except for specified.
Finally, $L$ controls the order of the regularized Dirichlet and 
Lagrange kernels and is set to 50 in calculations. 
It should be point out 
that, the selection of $L ~(L \geq M)$ is 
independent of the grid used in the present computation.

Three different problems in computational electromagnetics, 
waveguide analysis, electrostatic analysis and 3D
electromagnetic wave propagation, 
are selected to test the DSC algorithm. 
In all cases, we utilize the DSC algorithm to solve 
differential equations.
Details of these computations are 
described in the next three subsections. 
Double precision is used in all calculations.

\subsection{Waveguide analysis}

{\sc Eigenmode analysis}~~
The propagation of uniform plane waves is characterized by 
the transverse electromagnetic (TEM) nature of the wave.
The electric and magnetic field intensities are orthogonal to 
each other and to the direction of wave propagation. 
In a waveguide, the wave propagation is 
guided along given directions with particular characteristics
determined by the structure of the wave guide. 
Although it is assumed that boundaries are perfectly
conducting, the material within the waveguide is 
arbitrary and may include lossy or perfect dielectrics, 
conductor, etc.
The properties of a waveguide with complex geometries 
can only be simulated numerically. To test the DSC algorithm 
for CEM applications, we first consider the problem of finding 
the eigenvalues that determine
the parameters of waveguide modes, resonant frequencies of 
resonators, and many other physical parameters. 
Consider a rectangular waveguide, with propagation in the $z$-direction.
For the transverse magnetic (TM) mode, the four field components 
$E_x, E_y, H_x$ and $H_y$ can be expressed in terms of $E_z$. 
In turn, $E_z$ can be written as
%\vspace*{-.2cm}
\begin{equation}\label{waveg1}
E_z(x,y,z,t)=E(x,y)e^{i(\omega t-\alpha z)},
\end{equation}
%\vspace*{-.4cm} 
where $E$ satisfies 
%\vspace*{-.2cm}
\begin{equation}\label{waveg2}
{\partial^2 E\over \partial x^2}
+
{\partial^2 E\over \partial y^2}+k^2E=0,
\end{equation} 
%\vspace*{-.5cm}
and vanishes on the boundary ($E=0$) for the TM modes. 
The eigenvalue $k^2$ determines the phase parameter $\alpha$
through
%\vspace*{-.2cm}
\begin{equation}\label{waveg4}
k^2=\omega^2\nu\epsilon-\alpha^2,
\end{equation} 
%\vspace*{-.5cm}
where $\epsilon$ and $\nu$ are dielectric constant and magnetic 
permeability, respectively.

To simplify the problem further, 
we take the computational domain as $[0,10\pi]\times[0,10\pi]$.
Three DSC kernels are used to compute the eigenvalue problem
at three different grid sizes ($36^2, 24^2$ and $12^2$) so that 
the rate of convergence of the DSC approach for waveguide 
analysis can be examined. The computational bandwidth is set to 
the same as the grid size in this problem.
When bandwidthes are $M=36, 24$ and $12$, 
the values of $\sigma/\Delta$ are chosen as 4.2,  3.2  and 2.65
for the RSK and RDK and 2.95,  2.75 and 2.3 for the RLK. 
Results are obtained by a 
direct matrix diagonalization using
a standard eigenvalue solver.
Selected eigenvalues and absolute errors of the present 
numerical calculation are listed in TABLE I.

We first note that the RSK and RDK have very similar 
performance in all computations. Hence, we only 
need to consider  one of them in the rest of the 
computation. However, the performance of the RLK differs  
slightly from that of the RSK.

For $N=N_x=N_y=12$, there are only a maximum of  2.4 grid
points per wavelength in each dimension. 
Obviously, if this is the case, it is impossible for any 
of the three kernels or any other method to catch up with 
the fast oscillations in high-order 
eigenmodes.  All three kernels deliver excellent
results when the mesh is refined to $24^2$ (about 
4.8 grid points for each wavelength in each dimension).
Machine precision is attained for the first 100 
eigenmodes when the mesh is further refined to $36^2$.

{\sc Field distribution}~~
In practical applications, the geometry of a waveguide 
can be very complex so as to achieve certain required 
properties. Therefore, it is important for a 
computational algorithm to be flexible enough in handling 
complex boundary conditions. In this study, we 
consider two waveguides with geometries having 
a T-shape and an E-shape.
These boundaries are designed to test the capability 
of the DSC algorithm 
for complex waveguide simulation. We are 
interested the field distribution of 
TM (and/or TE) modes of these waveguides. 
The RSK is used for both cases and the value  of 
$\sigma/\Delta$ is set to 4.8. A total of $50$ 
DSC-collocation points are used in each dimension. 
Field distributions across the cross section for 
eigenmodes 1, 2, 3 and 4 are shown in Fig. 1a. 
These eigenmodes are quite localized and do not 
give the full shape of waveguide confining geometry. 
However, as plotted in Fig. 1b, higher-order eigenmodes, 
17, 18, 19 and 20, clearly reflect the waveguide geometry. 
Similar features are observed in the E-shape waveguide
as shown in Fig. 2. The first four modes of  E-shape 
waveguide are given in Fig. 2a.  Modes 17, 18, 19 
and 20  are depicted in Fig. 2b.

\subsection{Electrostatics analysis}

Potential function provides a useful alternative representation 
of electromagnetic fields. In many situations, it is more 
convenient to use an auxiliary function to describe  
field properties rather than to use full electric and 
magnetic field variables. In particular, a scalar potential can be 
used in a source-free domain and/or in domains which contain
scalar sources such as charge densities.

The electric scalar potential is one of the most widely 
used potential functions. If the rate of change of the 
magnetic flux vanishes, the electric field is irrotational
\begin{equation}
\nabla\times {\bf E}=0.
\end{equation}
In such a case, the  electric field can be written as the 
gradient of a scalar potential 
\begin{equation}\label{pot1}
{\bf E}=-\nabla V.
\end{equation}
For static fields, the divergence of Eq. (\ref{pot1})  
gives
\begin{equation}\label{pot2}
\nabla^2 V=-{\rho\over\varepsilon},
\end{equation}
where $\rho$ is the volume electric charge density
and $\epsilon$ is the electric permittivity. Here,
the Maxwell equation 
$\nabla\cdot{\bf E} ={\rho\over\varepsilon}$
is used for the derivation. Equation (\ref{pot2}) can be used as a 
starting point for many useful theoretical descriptions.

To demonstrate the DSC algorithm for electrostatic calculations,
we first consider a simple problem in two spatial dimensions.
The geometry of the problem 
 consists of a conducting square box of 
$1m\times 1m$ in the $x-y$ plane. The box is infinitely 
long in the $z-direction$ and thus, this can be regarded
as a two-dimensional problem.
The boundary condition of the problem is that the top side 
is isolated from all other sides with a
potential of $V(x,1)=10V$. All other three
sides are connected to the zero potential 
($V(0,y)=V(1,y)=V(x,0)=0$).  The material inside 
the box is the free space with permittivity 
$\varepsilon_0={1\over36\pi}\times 10^{-9}$.

We next consider an extension to the first problem. 
The same box with the same  boundary 
condition is used except that there is a $0.18m\times0.16m$ area of 
charge density inside  the box.   The charge density 
is uniformly distributed 
within a material of relative permittivity 
$\varepsilon_r=100$ and has a value of 
$1.0\times10^{-7} [C/m^2]$.

In both cases, we choose 32 points ($N=32$) in each dimension 
($\Delta x=\Delta y=0.0322m$). The DSC-collocation 
method is very convenient for handling these potential 
function problems.
The solution is obtained by solving the coupled algebraic 
equations (implicit scheme) with all boundary conditions 
being implemented in the operator part $({\cal L})$ of 
the coupled equations. The charge density is treated as a 
source term $(f)$ in the coupled equations.
For the first case, there is no charge density in the 
inner region and  the Laplace equation is directly solved. 
The results are plotted in Fig. 3.

For the second case, Eq. (\ref{pot2}) is modified as 
\begin{equation}\label{pot3}
\nabla^2 V(x,y) =-{\rho\over\varepsilon_0\varepsilon_r}.
\end{equation}
The charge density at the area of $x: 0.41m-0.59m,~ 
y: 0.72m-0.88m$ 
and the boundary conditions are directly imposed 
in a set of coupled algebraic equations, which are 
generated by using the DSC-collocation formulation. 
A standard linear algebraic equation solver is used 
to attain results as plotted in Fig. 4.

\subsection{Electromagnetic wave propagation}

In an uncharged homogeneous media, propagation of  
electromagnetic waves is governed by 
the wave equation
%\vspace*{-.2cm}
\begin{equation}\label{wave1}
{\partial^2 W\over\partial t^2}
-{1\over \epsilon\nu}\nabla^2W=0,
\end{equation}
%\vspace*{-.5cm}
where $W({\bf r},t)$ can be either the electric field or 
the magnetic field. To illustrate the potential of the DSC 
algorithm for electromagnetic wave simulations, 
we consider the following initial value problem
%\vspace*{-.2cm}
\begin{equation}\label{wave4}
W(x,y,z,0)=\sin(\alpha_xx+\alpha_yy+\alpha_zz),
\end{equation}
with periodic boundary conditions in all directions.
This problem has an  analytical solution given by
\begin{equation}\label{wave5}
W(x,y,z,t)=\sin\left[\alpha_xx+\alpha_yy+\alpha_zz+
\sqrt{\alpha_x^2 +\alpha_y^2+ \alpha_z^2\over\epsilon\nu}t\right].
\end{equation}
To further simplify notations, we set 
$\alpha_x= \alpha_y = \alpha_z=\epsilon\nu=1$.
The initial wave is propagated in a cubic size
$(10\pi)^3$. 
The fourth order explicit Runge-Kutta scheme is used for
time discretization. Sufficiently small
time increments ($\Delta t$) are used so that the 
major errors are due to 
the spatial discretization. The RSK and RLK
are used in the finite difference manner, i.e. 
according to Eq. (\ref{Eq10}), 
and DSC weights are calculated by using Eq. (\ref{Eq11})
once for the whole computation. 
We compute the problem by using a variety of grid points $N$ 
($N=N_x=N_y=N_z$) and computational bandwidth $M$ ($M=M_x=M_y=M_z$) 
values to test the accuracy and rate of convergence of the DSC 
algorithm and to demonstrate the unified features of global and 
local methods in the DSC algorithm. When $M=$36, 24, 12, and 6
the values of $\sigma/\Delta$ 
are 4.2,  3.2, 2.65 and 2.0 for the RSK 
and 2.95  2.75, 2.3 and 1.85 for the RLK. 
The $L_\infty$ errors at a number of time units
are listed in TABLE II.

When $N=24$, the number of grid points per wavelength 
is 4.8 in each dimension (Note that typically the 
Yee algorithm uses about 18 points per wavelength). 
This is a typical case
of under sampling and it is very difficult to achieve 
high computational accuracy by any means. We choose
three different $M$ values. For $M=6$, results 
from the two DSC kernels are accurate enough. 
The accuracy improves as the $M$ value is increased from  6 to
12 and 24. The last case, $M=N=24$, corresponds to a global
treatment (i.e., a ``global finite difference'') 
and its results are significantly better.

In the case of $N=36$, it is on an average about 
7.2 grid points per wavelength in each dimension. 
The computational accuracy increases as the 
computational bandwidth increases. It is noted that
at $M=12$ and $M=24$, the accuracy does not improve 
much from the previous computations ($N=24$; $M=12, 24$)
and is obviously limited by the computational bandwidth.
However, it is interesting to note that the DSC algorithm 
attains extremely high accuracy at the global limit of
$M=N=36$. To the best of our knowledge,
these are the highest precision results obtained for 
this problem, so far.

\section{Conclusion}

In conclusion, this paper introduces the
discrete singular convolution (DSC)
algorithm for computational electromagnetics (CEM).
The computational philosophy of the DSC algorithm 
is studied. Many sequences of approximations to 
the delta distribution, the ``universal 
reproducing kernel'',  are constructed either as 
bandlimited reproducing kernels or as approximate 
reproducing kernels. A regularization procedure
based on the distribution theory is utilized to 
improve the regularity and localization of the DSC 
kernels. Systematic treatments of 
 derivatives and boundary conditions are proposed, 
as they are required for a computational algorithm.

We explore the unified features of DSC algorithm 
for solving differential equations in the framework 
of the method of weighted residuals. It is found that several
computational methods, such as global, local,
Galerkin, collocation and finite difference methods, can be 
derived from a single starting point by using 
DSC trial functions and test functions. 
The unification of local and global methods is realized via 
the DSC approach. It is well known that 
accuracy is crucial to many scientific 
computations, such as  simulation of turbulence
and high-frequency analysis of radar cross-section.
Whereas,  small-matrix-band approximations are vital to 
large scale computations. The DSC 
 approach provides a controllable numerical accuracy 
by an appropriate selection of the matrix bandwidth. 
Therefore, by using the DSC algorithm, 
the computational efficiency can be easily optimized
against accuracy, convergence, stability and simplicity.
A collocation algorithm is deduced from a Galerkin
statement and is called a Galerkin-induced collocation. 
The equivalence of Galerkin and 
collocation enables us to  evaluate 
these two conventional numerical algorithms on an equal
footing. The connection is made between finite 
difference  and other methods, such as
collocation and Galerkin. The DSC algorithm can be regarded
as a generalized finite difference method, which becomes a 
``global finite difference method'' by an appropriate  
choice of the computational bandwidth.
These results can be used as a guide 
for the selection of computational methods and for the
design of numerical solvers for practical applications.

The potential of the DSC algorithm for computational 
electromagnetics is explored by performing  a number of 
numerical experiments. 
The first example is about waveguide eigenmode analysis,
an eigenvalue problem in a square domain of $(10\pi)^2$.
This problem has an analytical solution and hence 
can be used for testing the accuracy of potential 
numerical methods. We formulate the problem in 
the DSC-collocation approach
and attain results by using 
three sets of grid points, $12^2, 24^2$ and $36^2$, to 
examine the speed of convergence of the DSC algorithm. 
Three typical DSC kernels, a regularized Shannon's kernel,
a regularized Dirichlet kernel, and a regularized Lagrange 
kernel are utilized to illustrate the usefulness of 
the proposed algorithm. It is found that 
the accuracy improves dramatically as the number of grid 
point is increased. In particular, it reaches the machine
precision for the first 100 eigenmodes when the mesh is
refined to $36^2$. Having built enough confidence for eigenmode 
computations, we have tested the ability of the DSC algorithm 
in handling complex boundary conditions. 
In this regard, we have used the DSC algorithm for TM field 
distributions of T-shape and E-shape waveguides.
We found that the first few eigenmodes are very 
localized due to complexity in geometry. However,
higher-order eigenmodes exhibit global characteristics.
The confining geometry of the waveguide can be seen from 
mode patterns.

In the second example, we consider potential function 
analysis of electrostatics in two spatial dimensions.
Computationally, this is a non-uniform boundary value
problem. Only one side of the four boundaries has a 
non-zero potential value. We have treated the problem 
via the DSC-collocation approach with a mesh of $32^2$.
Since the DSC algorithm is a local method in general,
boundary conditions and charge density are easily 
implemented in a set of coupled algebraic equations. 
The DSC approach has proven to be very efficient and 
robust for this problem.

In the last example, we have tested the DSC algorithm for
three-dimensional electromagnetic wave propagations. 
The wave equation is integrated over a long time 
with periodic boundary conditions. Since the problem 
admits an analytical traveling wave solution, the 
performance of the DSC algorithm can be examined objectively. 
This problem is rather sensitive to the 
computational accuracy and stability. The temporal
discretization is carried out by using the standard
fourth order explicit Runge-Kutta scheme. The DSC 
algorithm is utilized for spatial discretization
in the domain of $(10\pi)^3$.
The unified feature of the DSC algorithm for both 
global and local approximations is illustrated in
this problem by varying the computational bandwidth
$M$ for a given matrix size $N$. In this calculation, 
the DSC algorithm is used as a generalized finite 
difference method. In particular, when the bandwidth equals 
the matrix length ($M=N$), the DSC algorithm gives rise to 
an interesting ``global finite difference'' treatment.
The DSC algorithm can provide controllable 
accuracy by varying the parameter $M$, hence,
one of the advantages of the DSC algorithm is its 
robustness, and allows the selection of 
desirable accuracy for a given problem without any need to
change one's computer code.
All numerical solutions are
found to be stable over a long time integration. A sharp 
improvement in the numerical accuracy is observed when the mesh 
is refined from $24^3$ to $36^3$, i.e. from 4.8 points 
per wavelength to 7.2 points per wavelength 
in each dimension. It is found that 
the accuracy of the DSC algorithm is at least 12 
significant figures up to 22 time units
when the computational bandwidth is relatively high
($M=N=36$). This result is consistent with our previous
theoretical error estimation\cite{qian99}. 
Present work indicates that the DSC algorithm is 
a promising and potential approach for 
computational electromagnetics. 

\vspace*{.2cm}

\centerline{\bf Acknowledgment}
 
{This work was supported in part by the National University 
of Singapore and in part by the National Science and Technology 
Board of the Republic of Singapore.}

%\end{document}       
 
\bibliographystyle{plain}

\newpage

\centerline{\bf Figure Caption}

FIG. 1a. Contour and mesh plots of 
the first four eigenmodes of the  T-shape waveguide.

FIG. 1b. Contour and mesh plots of 
eigenmodes 17, 18, 19 and 20 of the  T-shape waveguide.

FIG. 2a. Contour and mesh plots of 
the first four eigenmodes of the  E-shape waveguide.

FIG. 2b. Contour and mesh plots of 
eigenmodes 17, 18, 19 and 20 of the E-shape waveguide.

FIG. 3. Contour and mesh plots  of the potential field distribution
  in a square box with non-uniform boundaries.  

FIG. 4. Contour and mesh plots  of potential field 
       distribution   in a square box with non-uniform boundaries
      and  an area of charge density.

\newpage

\begin{center}

\begin{table}
\caption{Errors for waveguide modes}
\begin{center}
\begin{tabular}{c|c|c|c|c} 
$N$& Eigenmode&      ~ RSK &    RDK      &   RLK \\ \hline
12 &   10   &    8.4(-8)   &   8.8(-8)   &   2.1(-7)\\
   &   20   &    5.2(-7)   &   5.1(-7)   &   4.7(-7)\\
   &   30   &    1.9(-2)   &   1.9(-2)   &   2.0(-2)\\
   &   40   &    1.1(-1)   &   1.1(-1)   &   1.1(-1)\\
   &   50   &    1.5(-1)   &   1.5(-1)   &   1.5(-1)\\
   &   60   &    1.9(-1)   &   1.9(-1)   &   1.9(-1)\\
   &   70   &    2.7(-1)   &   2.7(-1)   &   2.7(-1)\\
   &   80   &    3.2(-1)   &   3.2(-1)   &   3.2(-1)\\
   &   90   &    3.7(-1)   &   3.7(-1)   &   3.7(-1)\\
   &  100   &    4.2(-1)   &   4.2(-1)   &   4.2(-1)\\
24 &   10   &    1.8(-14)  &   3.1(-15)  &   9.7(-13)\\
   &   20   &    9.5(-15)  &   1.6(-14)  &   3.2(-12)\\
   &   30   &    8.5(-13)  &   8.2(-13)  &   1.4(-12)\\
   &   40   &    8.7(-13)  &   8.5(-13)  &   5.1(-13)\\
   &   50   &    3.8(-10)  &   3.7(-10)  &   1.3(-11)\\
   &   60   &    1.4(-8)   &   1.3(-8)   &   3.5(-10)\\
   &   70   &    8.0(-8)   &   7.8(-8)   &   6.3(-9) \\
   &   80   &    8.0(-8)   &   7.8(-8)   &   6.3(-9) \\
   &   90   &    7.6(-10)  &   7.5(-10)  &   3.0(-11)\\
   &  100   &    6.4(-9)   &   6.3(-9)   &   3.6(-10)\\
36 &   10   &    1.2(-14)  &   2.1(-14)  &   2.9(-14)\\
   &   20   &    4.3(-14)  &   2.1(-14)  &   1.2(-15)\\
   &   30   &    1.4(-14)  &   2.8(-14)  &   7.7(-15)\\
   &   40   &    2.2(-16)  &   3.9(-14)  &   5.0(-15)\\
   &   50   &    2.0(-14)  &   2.4(-14)  &   9.0(-14)\\
   &   60   &    7.4(-14)  &   2.0(-14)  &   3.8(-14)\\
   &   70   &    2.4(-14)  &   7.3(-15)  &   2.9(-14)\\
   &   80   &    2.2(-14)  &   4.2(-14)  &   6.7(-15)\\
   &   90   &    1.3(-14)  &   6.4(-15)  &   4.4(-14)\\
   &  100   &    3.7(-14)  &   1.8(-15)  &   3.3(-14)\\
\end{tabular}
\end{center}
\end{table} 

\end{center}

\newpage

\begin{center}

\begin{table}
\caption{ $L_\infty$ errors for solving the 3D wave  
            equation}
\begin{center}
\begin{tabular}{c|c|c|c|c|c} 
$N$ & $M$ & $\Delta t$ & $t$ & RSK & RLK  \\ \hline
24 & 6  & 0.5 &  1      &   6.2(-4)  & 5.1(-4)\\
   &    &     &  4      &   1.1(-3)  & 1.5(-3)\\
   &    &     &  7      &   1.7(-3)  & 2.4(-3)\\
   &    &     &  10     &   2.3(-3)  & 3.3(-3)\\
   &    &     &  13     &   2.8(-3)  & 4.3(-3)\\
   &    &     &  16     &   3.5(-3)  & 5.2(-3)\\
   &    &     &  19     &   4.1(-3)  & 6.1(-3)\\
   &    &     &  22     &   4.8(-3)  & 7.0(-3)\\
   & 12 & 0.2 &  1      &   1.2(-5)  & 1.2(-5)\\
   &    &     &  4      &   2.8(-5)  & 3.1(-5)\\
   &    &     &  7      &   4.5(-5)  & 5.1(-5)\\
   &    &     &  10     &   6.1(-5)  & 7.1(-5)\\
   &    &     &  13     &   7.8(-5)  & 9.0(-5)\\
   &    &     &  16     &   9.4(-5)  & 1.1(-4)\\
   &    &     &  19     &   1.1(-4)  & 1.3(-4)\\
   &    &     &  22     &   1.3(-4)  & 1.5(-4)\\
   & 24 & 0.05&  1      &   4.3(-8)  & 5.0(-8)\\
   &    &     &  4      &   7.6(-8)  & 2.0(-7)\\
   &    &     &  7      &   1.1(-7)  & 3.4(-7)\\
   &    &     &  10     &   1.4(-7)  & 4.9(-7)\\
   &    &     &  13     &   1.8(-7)  & 6.4(-7)\\
   &    &     &  16     &   2.1(-7)  & 7.8(-7)\\
   &    &     &  19     &   2.4(-7)  & 9.3(-7)\\
   &    &     &  22     &   2.7(-7)  & 1.1(-6)\\
36 & 12 & 0.2 &  1      &   1.3(-5)  & 1.3(-5)\\
   &    &     &  4      &   4.8(-5)  & 4.7(-5)\\
   &    &     &  7      &   8.6(-5)  & 8.4(-5)\\
   &    &     &  10     &   1.2(-4)  & 1.2(-4)\\
   &    &     &  13     &   1.6(-4)  & 1.5(-4)\\
   &    &     &  16     &   2.0(-4)  & 1.9(-4)\\
   &    &     &  19     &   2.3(-4)  & 2.2(-4)\\
   &    &     &  22     &   2.6(-4)  & 2.5(-4)\\             
   & 24 & 0.05&  1      &   5.2(-8)  & 5.2(-8)\\
   &    &     &  4      &   2.1(-7)  & 2.1(-7)\\
   &    &     &  7      &   3.6(-7)  & 3.6(-7)\\
   &    &     &  10     &   5.2(-7)  & 5.2(-7)\\
   &    &     &  13     &   6.8(-7)  & 6.8(-7)\\
   &    &     &  16     &   8.2(-7)  & 8.2(-7)\\
   &    &     &  19     &   9.8(-7)  & 9.8(-7)\\
   &    &     &  22     &   1.1(-6)  & 1.1(-6)\\
   & 36 &0.001&  1      &   1.3(-14) & 1.4(-14)\\
   &    &     &  4      &   3.0(-13) & 3.0(-13)\\
   &    &     &  7      &   7.3(-13) & 7.3(-13)\\
   &    &     &  10     &   3.5(-14) & 2.7(-14)\\
   &    &     &  13     &   1.7(-12) & 1.7(-12)\\
   &    &     &  16     &   3.3(-12) & 3.3(-12)\\
   &    &     &  19     &   3.9(-13) & 4.1(-13)\\
   &    &     &  22     &   4.1(-12) & 4.1(-12)\\
\end{tabular}
\end{center}               
\end{table} 

\end{center}

\end{document}